\documentclass{amsart}
\usepackage[foot]{amsaddr}
\usepackage{lmodern}
\usepackage[utf8]{inputenc}
\usepackage[T1]{fontenc}
\usepackage{graphicx}
\usepackage{subcaption}
\usepackage{amsfonts}
\usepackage{amsmath}
\usepackage{amssymb}
\usepackage{float}
\usepackage[dvipsnames]{xcolor}
\usepackage[colorlinks=true,allcolors=PineGreen]{hyperref}
\usepackage{url}
\usepackage{doi}
\usepackage{orcidlink}
\usepackage[normalem]{ulem}
\usepackage{comment}

\providecommand{\llabel}[1]{}
\providecommand{\lref}[1]{}

\usepackage[section]{placeins}

\usepackage[backend=biber, style=numeric, giveninits=true, uniquename=init, uniquelist=false, sorting=nyt]{biblatex}

\DeclareSourcemap{
  \maps[datatype=bibtex]{
    \map[overwrite]{
      \step[fieldsource=doi, final]
      \step[fieldset=url, null]
      \step[fieldset=eprint, null]
    }  
  }
}

\AtEveryBibitem{\clearfield{issn}}
\AtEveryCitekey{\clearfield{issn}}
\AtEveryBibitem{\clearfield{isbn}}
\AtEveryCitekey{\clearfield{isbn}}

\DeclareMathOperator{\inv}{Inv}
\DeclareMathOperator{\cl}{cl}
\DeclareMathOperator{\interior}{int}

\theoremstyle{definition}
\newtheorem{defin}{Definition}

\title[Analysis of Global Dynamics in the Andrecut--Kauffman Model]{Topological--numerical analysis of global dynamics in the discrete-time two-gene Andrecut--Kauffman model}

\author{Dorian Falęcki $^{1,\dag,\orcidlink{0009-0006-3451-9128}}$}
\address{$^1$ Faculty of Applied Physics and Mathematics,
Gda\'{n}sk University of Technology,
ul.~Gabriela Narutowicza 11/12, 80-233 Gda\'{n}sk, Poland}
\address{$^\dag$ Both authors contributed equally to the research}

\author{Mikołaj Rosman $^{1,\dag,\orcidlink{0009-0002-7760-2491}}$}

\author{Michał Palczewski $^{2,\orcidlink{0009-0007-3775-9589}}$}
\address{$^2$ Doctoral School and Faculty of Applied Physics and Mathematics,
Gda\'{n}sk University of Technology,
ul.~Gabriela Narutowicza 11/12, 80-233 Gda\'{n}sk, Poland}

\author{Paweł Pilarczyk $^{3,\orcidlink{0000-0003-0597-697X}}$}
\address{$^3$ Faculty of Applied Physics and Mathematics \& Digital Technologies Center,
Gda\'{n}sk University of Technology,
ul.~Gabriela Narutowicza 11/12, 80-233 Gda\'{n}sk, Poland}

\author{Agnieszka Bartłomiejczyk $^{4,\ast,\orcidlink{0000-0003-1086-8631}}$}
\address{$^4$ Faculty of Applied Physics and Mathematics \& BioTechMed Center,
Gda\'{n}sk University of Technology,
ul.~Gabriela Narutowicza 11/12, 80-233 Gda\'{n}sk, Poland}
\address{$^\ast$ Corresponding author}
\thanks{This is the authors' peer reviewed, accepted manuscript (AAM; postprint)
of the article published on July 29, 2026 in
\textit{Chaos: An Interdisciplinary Journal of Nonlinear Science},
Vol.~36, Issue~7, 073150; \doi{10.1063/5.0333093}.
This version of the article may be downloaded for personal use only.
Any other use requires prior permission of the authors and AIP Publishing.}

\begin{document}

\begin{abstract}
We conduct a topological--numerical analysis of global dynamics in a discrete-time two-gene Andrecut--Kauffman model. This model describes gene expression regulation through nonlinear interactions.
We use a numerical method to construct Morse decomposition of the system across a wide range of parameters at a fixed finite resolution in the state space and in the parameter space, both spaces split into uniform rectangular grids (a technique also called ``pixelation'').
We obtain qualitative results by effectively computing the Conley indices of the constructed isolating neighborhoods that form the Morse decomposition.
We represent the Morse decomposition and connecting orbits by a directed acyclic graph.
We introduce pictograms to convey the information provided by the Conley index in an easy to understand schematic way.
We show and analyze bifurcations captured using this technique.
We call this method CMAD for short (Conley-Morse graphs for the Analysis of Dynamics).
The main advantage of our method is that it finds isolating neighborhoods of both stable and unstable invariant sets and that it provides validated (rigorous) numerical results: we actually obtain computer assisted proof that the constructed sets are indeed isolating neighborhoods of Morse sets in a certain Morse decomposition of the system. In particular, this means that we have captured all the interesting dynamics within the analyzed range of the phase space perceived at the given finite resolution.
We also conduct numerical simulations aimed at showing the location of attractors in the isolating neighborhoods found.
The results demonstrate the usefulness of topological methods in understanding the global structure of dynamics at finite (coarse) resolution in an applied dynamical system depending on a few parameters, like the gene regulatory model that we analyze.
\end{abstract}



\maketitle


\textbf{Simplified mathematical models of complex processes like gene expression regulation help us to gain insight into some aspects of the sophisticated phenomena that cannot be directly modeled in full detail. However, dynamics generated by such models, like the discrete-time two-dimensional Andrecut--Kauffman model, may still be too difficult to fully understand using the analytic approach. Therefore, numerical simulations are often applied instead. Unfortunately, such methods are typically limited to finding stable dynamics only and are subject to inaccuracies (numerical errors) that undermine their reliability. We show an application of a specific topological--numerical method to obtain classification of the dynamics found in the Andrecut--Kauffman model for a wide range of parameters. Our approach uses rigorous numerical methods and topological tools such as the Conley index and Morse decompositions. We obtain a picture of the dynamics that encompasses all stable and unstable invariant sets. The results are obtained at a relatively coarse level (limited to an apriori chosen finite resolution), but they are mathematically reliable\llabel{reliable} in the sense that we obtain a computer assisted proof that the constructed sets are indeed isolating neighborhoods of a true Morse decomposition of the system. As a consequence, we know that we have not missed any invariant sets in the phase space, stable or unstable, including attractors with very small attraction basins that can be missed in simple numerical simulations. We illustrate and discuss the results obtained and show how these results may be refined by means of classical numerical simulations to obtain a comprehensive and refined picture of the dynamics and bifurcations.}

\section{Introduction}
\label{sec:intro}

Understanding the dynamics of nonlinear discrete systems is a fundamental problem in contemporary applied mathematics and computational science. Such systems may exhibit a broad spectrum of behaviors, including stable and unstable fixed or periodic points, bistability, bifurcations, and deterministic chaos.

Classical analytic techniques usually focus on local properties of the dynamics, such as fixed or periodic points and their linear stability. Although these methods are effective in many settings, they often fail to provide a complete description of the global organization of trajectories in phase space, especially in systems with
unstable invariant sets.

To overcome these limitations, increasing attention has been given to the topological and rigorous computational methods for nonlinear dynamics. Tools such as the Conley index, Morse decompositions, graph representations, and computational homology allow one to detect stable and unstable invariant sets, describe connections between them, and identify structural changes in the dynamics. In contrast to standard numerical simulations alone, these approaches provide mathematically validated information about the global structure of the system:
they provide results in the form of a computer-assisted proof.
The effectiveness of early methods using this approach was shown in the 1990s and included the computer-assisted proof of the existence of chaotic dynamics in the Lorenz equations~\cite{MM1995chaos}, recognized by the \emph{Encyclopedia Britannica} as one of the four most important achievements in mathematics in 1995. These methods have been undergoing continuous development since then; see e.g. \cite{Dey2022,Kapela2021,Mischaikow2016} for some more recent directions.

In this paper, we consider the discrete two-gene Andrecut–Kauffman model \cite{andrecut_main}, derived from the classical Boolean networks introduced by Kauffman in 1969 \cite{KAUFFMAN1969437}. Despite its simple threshold-based formulation, the model generates a rich variety of dynamical phenomena ranging from stable steady states and oscillations to bistability and deterministic chaos \cite{Andrecut_2005, Andrecut_2006}.

Previous studies of this model have focused mainly on local dynamical features or numerical exploration of selected parameter regions. However, a comprehensive understanding of the system requires methods capable of describing the full phase-space organization and transitions between dynamical regimes. The approach proposed here enables us to capture relationships among trajectories and their topological structure in a systematic way.

In this paper, we apply a combined topological and numerical framework to analyze the global dynamics of the discrete two-gene Andrecut–Kauffman model. We employ discrete graph algorithms and computational homology tools to identify and classify stable and unstable invariant sets, bifurcations, and structural transitions. Using rigorous numerical methods, the obtained results are mathematically reliable\llabel{reliableMeaning} in the sense that they constitute computer-assisted proofs of the reported features. The computations are conducted at limited resolution using the idea of pixelation of the phase space and the parameter space, explained in Sec.~\ref{sec:numMorse}. Therefore, the information we obtain on the location of sets and bifurcations is limited to the chosen resolution. Therefore, we use complementary numerical simulations to provide additional insight into the geometry of attractors (particularly in chaotic parameter regimes) and more precise location of bifurcations that our method detects.

An important reference point for the present work is our previous study \cite{rosman-2025}, where bistability and chaos in the same model were investigated primarily through numerical simulations over a broad parameter range. In the current paper, we extend those results by incorporating computational topology, thereby obtaining a more complete and global description of the dynamics that also includes unstable invariant sets.

The results presented here demonstrate the effectiveness of combining topological methods with numerical simulations in the comprehensive characterization of discrete systems. More broadly, they illustrate how rigorous computational tools can complement classical dynamical systems' analysis in the study of nonlinear models.

The paper is organized as follows. We introduce the discrete two-gene Andrecut–Kauffman model in Sec.~\ref{sec:model}, review the current state of the art in Sec.~\ref{sec:state}, and define the scope and explain the originality of the presented contribution in Sec.~\ref{sec:contrib}. In Sec.~\ref{sec:method}, we focus on the methodology, introduce the topological--numerical tools used to analyze the dynamics of the system, and briefly introduce their theoretical foundations, with particular emphasis on the Conley index and Morse decompositions. In Sec.~\ref{sec:results}, we show the main results and explain their interpretation. We describe the chosen set of parameters, the numerical procedures employed, and the resulting dynamical structures found in the system. We also examine the influence of parameter values and phase space resolution on the observed dynamics. The paper concludes with a summary of the main findings and final remarks.

\subsection{The discrete-time two-gene Andrecut--Kauffman model}
\label{sec:model}

We consider a simplified system consisting of two genes in which the complex processes of transcription and translation are represented by a single effective reaction. Although this approximation omits many biochemical details, it retains the essential mechanisms necessary for the analysis of gene expression and its regulation.
The expression of each gene is modeled by a reaction in which RNA polymerase and the gene's promoter lead to the production of a protein monomer. Regulatory control is incorporated by allowing multimerized forms of the proteins to bind to promoter regions, thereby modulating expression levels. The model also accounts for the multimerization process itself, in which protein monomers combine to form functional transcription factors, and includes degradation reactions that remove protein monomers from the system over time.
These biochemical interactions are described by the following pair of difference equations that induce a discrete-time dynamical system\llabel{semiDynam} (note that the map is not invertible):
\begin{equation}
\tag{AK}
\label{eq:model}
\left\{\;
\begin{aligned}
x_{t+1} &= \frac{\alpha_1}{1+(1-\varepsilon)x^n_t+\varepsilon y^n_t} +\beta_1 x_t,\\
y_{t+1} &= \frac{\alpha_2}{1+\varepsilon x^n_t+(1-\varepsilon) y^n_t}+\beta_2 y_t.
\end{aligned}
\right.
\end{equation}
Here, $x_t$ and $y_t$ denote the concentrations of the two proteins at the time $t \in \mathbb{N}$.
Protein concentrations vary in time according to the balance between their production, modulated by regulatory interactions, and their degradation. This system depends on six\llabel{key_removed} parameters: the expression strength $\alpha_1$ and $\alpha_2$ of the two genes, the degradation parameters $1-\beta_1$ and $1-\beta_2$, the degree of multimerization $n$, and the coupling parameter $\varepsilon\in [0,1]$.
In the literature, the dimension of the parameter\llabel{dimParamBegin} space of such models is often reduced by assuming symmetry in the expression and degradation parameters of both genes, that is, $\alpha_1 = \alpha_2$ and $\beta_1 = \beta_2$. Although the original formulation~\cite{andrecut_main} of the Andrecut--Kauffman model allows different values of these parameters for each gene, many authors adopted equal parameters as a deliberate mathematical simplification.\llabel{dimParamEnd}

\llabel{alpha12begin}In the present article, we investigate the behavior of the system for different values of $\alpha_1$ and $\alpha_2$ with the other parameters fixed, so that we work with a 2-dimensional parameter space. The possibility of $\alpha_1 \neq \alpha_2$ allows us to explore the impact of asymmetry in gene expression on discrete-time dynamics. Although our computational framework allows us to work with a higher-dimensional parameter space, the analysis and visualization of the results would be more challenging.\llabel{alpha12end}

\subsection{State of the art}
\label{sec:state}

In papers on the Andrecut--Kauffman model that have been published so far, it was shown that variations of the parameters lead to different dynamical behaviors, such as fixed points, periodic oscillations, chaos, hyperchaos, and quasiperiodicity. Let us provide a brief summary of some main results.

In \cite{sharma2019}, the authors analyze the model under the assumption that \(\beta = \beta_1 = \beta_2\). They generate bifurcation diagrams and attractor plots for \(n = 3, 4, 5\) with fixed values \(\alpha = 25\), \(\varepsilon = 0.1\). Their results show period-doubling cascades leading to chaos and eventual return to regularity. The range of values of \(\beta\) associated with complex dynamics expands with the degree of multimerization: for \(n = 3\), complex dynamics occurs for \(0.05 \leq \beta \leq 0.23\); for \(n = 4\), the range increases to \(0 \leq \beta \leq 0.22\); and for \(n = 5\), it is \(0 \leq \beta \leq 0.44\).\llabel{sharmaFolding}

In another study \cite{De_Souza2012-vs} that focuses on self-similarities in periodic structures, the parameters \(\beta_1\) and \(\beta_2\) are treated independently and the other parameters are fixed at $\alpha_1=\alpha_2=25$, $\varepsilon=0.1$, and $n=3$. The authors utilize bifurcation diagrams, Lyapunov exponent maps, and isoperiodic diagrams to characterize the parameter space with special focus on the structure of periodic windows. Numerical simulations reveal the presence of Arnold tongues and shrimp-shaped periodic regions; both structures are known to have the property of self-similarity. These periodic regions follow period-adding sequences, and some displaying relations are similar to the Fibonacci sequence. The ratios of successive periods converge toward the Golden ratio.

Further insights are provided in \cite{subramani-2023}, where bifurcation diagrams and maximum Lyapunov exponents are used to analyze the effects of independently varying \(\beta_1\) and \(\beta_2\). For fixed \(\beta_2 = 0.42\), varying \(\beta_1\) in the interval \([0, 0.25]\) reveals multiple disconnected regions of chaos, periodicity, and fixed points. Similarly, scanning \(\beta_2\) with fixed \(\beta_1 = 0.2\) uncovers alternating behaviors in \([0, 0.5]\). A similar analysis is performed for the parameters $\alpha$ and $\varepsilon$. The parameter $\alpha$ is considered in the range $[0, 150]$ for fixed values $n = 3$, $\varepsilon = 0.1$, $\beta_1 = 0.2$, and $\beta_2 = 0.42$. The values of $\alpha \in [0, 22.5]$ correspond to fixed-point behavior, reflecting stable steady-state dynamics. As $\alpha$ increases, chaotic behavior is observed over several disjoint intervals: $[22.5, 57]$, $[68.6, 84.2]$, $[91.8, 113.6]$, and $[128, 150]$, interwoven with periodic windows in the ranges $[57, 68.6]$, $[84.2, 91.8]$, and $[113.6, 128]$.  Similarly, the parameter $\varepsilon$ is varied in the interval $[0, 0.15]$, with other parameters  fixed at $n = 3$, $\alpha = 25$, $\beta_1 = 0.2$, and $\beta_2 = 0.42$. Chaotic behavior is found in the ranges $[0, 0.0136]$ and $[0.0686, 0.1356]$, while periodic behavior emerges in the intermediate interval $[0.0136, 0.0686]$. For higher values $\varepsilon \in [0.1356, 0.15]$, the system settles into fixed-point dynamics. This analysis further shows that these transitions are reflected in abrupt shifts in the largest Lyapunov exponent.
Parameter plane diagrams confirm the coexistence of multiple attractors and the complex shapes of their attraction basins.

\subsection{Our contribution}
\label{sec:contrib}

We apply a \llabel{modern_removed}method that combines rigorous numerics with computational topology to analyze the Andrecut--Kauffman model. This method was originally introduced in \cite{arai-2009} and was further developed, for example, in \cite{knipl-2015,pilarczyk-2023}. With this approach, we are able to expand upon the results from previous work by finding unstable invariant sets in addition to attractors. While finding stable sets is an easy task\llabel{easySim} to achieve by straightforward numerical simulations, unstable sets are virtually undetectable for such methods.

\begin{figure}[htbp]
    \centering
    \includegraphics[width=0.8\linewidth]{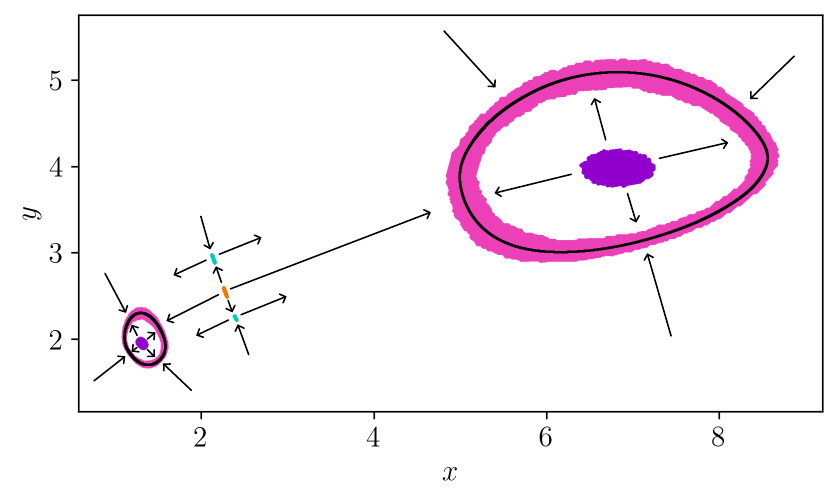}
    \caption{Phase space diagram for the Andrecut--Kauffman model \eqref{eq:model} with $\alpha_1=\alpha_2=25$, $\beta_1=0.18$, $\beta_2=0.42$, $\varepsilon=0.1$, and $n=3$, discussed in Sec.~\ref{sec:contrib}.
    }
    \label{fig:subramani}
\end{figure}

This difference is illustrated in the example shown in Fig.~\ref{fig:subramani}.
\llabel{twoCircBegin}In numerical simulations, orbits stabilize on a set that consists of two disjoint circles (shown in black). The orbit jumps between the circles and the second iterate of the map exhibits quasi-periodic motion on each circle.\llabel{twoCircEnd}
The authors of \cite{subramani-2023} were only able to find this stable invariant set in their numerical simulations.

\llabel{isolOnlyBegin}Instead of following individual trajectories, the method we apply works with sets with respect to an apriori declared resolution of the phase space (``pixelation''). It describes the dynamics by means of isolating neighborhoods that contain isolated invariant sets; see Sec.~\ref{sec:isol} for the definition of these terms. We therefore obtain coarse (but validated) information on the location of the invariant sets as well as possible connections between them, represented in the form of a Morse decomposition (see Sec.~\ref{sec:Morse}).
Isolating neighborhoods for all stable and unstable bounded invariant sets found with this method are shown in Fig.~\ref{fig:subramani} in different colors (note that some of them consist of two components). The sets in this figure have been drawn slightly thicker than their actual size for readability.
Arrows represent the stable and unstable directions of each set, as well as connecting orbits.

\llabel{focusBegin}Without additional computations, everything we know about the actual invariant sets comes from the properties of the constructed neighborhoods only, as explained in Sec.~\ref{sec:Conley}. In particular, if there is a stable fixed point, we do not know if it is a stable node or a stable focus. Such complementary information could be computed analytically or using classical numerical methods.\llabel{focusEnd}

Our approach has confirmed the existence and location of the attractor (its isolating neighborhood is shown in pink and comprises of two rings). It detected the topology of two circles (by computing their Betti numbers) and noticed the mapping between them (by computing the homomorphism on the homology groups). It additionally identified two period-$2$ orbits (their isolating neighborhoods are shown in purple and cyan\llabel{orange2purple})\llabel{isolOnlyEnd} and one smaller repelling fixed point (its neighborhood shown in red). The Conley index (see Sec.~\ref{sec:Conley}) provides information on the stability of all these sets (corresponding to the dimension of the unstable manifold), \llabel{arrowsBegin}indicated by the direction of the arrows. Our method also shows between which isolating neighborhoods connecting trajectories may exist, which is indicated in the figure by the arrows pointing from one set to another.\llabel{arrowsEnd} The computed information is rigorous (computer-assisted proof) and is in this sense different\llabel{different} from what can be obtained by plain numerical simulations. Such simulations\llabel{reliableArtifacts} might be vulnerable to rounding errors and yield numerical artifacts. Combining our method with numerical simulations provides a more complete overall understanding of the dynamics.

Using the results of our analysis, we propose partitioning the parameter space into subsets called \emph{continuation classes} (see Sec.~\ref{sec:continuation}) based on Morse decomposition of the phase space, which is significantly more granular than the currently used classification into chaotic and ordered regions. Detailed analysis of major continuation classes found for a chosen parameter regime illustrates what information can be extracted solely from Conley--Morse graphs and isolating neighborhoods. We complement this information with numerical approximations of attractors.

\section{The topological--numerical CMAD method}
\label{sec:method}

Let us provide an overview of the computational method that we apply and then recall some definitions and explain what information is provided by this method.
The method is based on Morse decompositions and the Conley index, and thus we shall call it CMAD\llabel{CMADname} as an abbreviation of ``Conley-Morse graphs for the Analysis of Dynamics.''
We refer to \cite{arai-2009,pilarczyk-2023} and references therein for additional background, explanations and details.

\subsection{Overview of the method}
\label{sec:overview}

The goal of the method is to understand the dynamics given by a difference equation in $\mathbb{R}^n$. The analysis is restricted to a bounded subset $B \subset \mathbb{R}^n$ of the phase space  perceived at a given resolution. We want to find the structure of dynamics in $B$ in the form of a finite number of sets at the given resolution containing all bounded invariant sets of the system (stable and unstable) and possible connecting orbits between them. The analysis is done over a specified bounded range of parameters, and in this way we capture bifurcations.

\llabel{missFineStructuresBegin}The analysis is conducted at a finite resolution using uniform grids in the phase space and in the parameter space; this technique can also be called ``pixelation.'' Therefore, we obtain a rather coarse view of the system and the information on the invariant sets is provided in the form of their neighborhoods and the stability of these neighborhoods. Although we miss the fine structures of dynamics \emph{inside} such neighborhoods and this method does not allow to find more precise\llabel{preciseLoc} location of the invariant sets than the rough bounds provided by the isolating neighborhoods, we obtain a complete picture of the dynamics at the chosen resolution in terms of these neighborhoods and possible connections between them, as well as bifurcations as the parameters vary. These neighborhoods can be used as a starting point for further (more thorough) analysis.\llabel{missFineStructuresEnd}

\llabel{rigNumBegin}The information we obtain is validated by means of rigorous numerics and algebraic topology (computer-assisted proof). This feature makes our results reliable in contrast to classical numerical simulations based on finite-precision floating-point arithmetic without rounding control which may lead to results different from what really happens in the actual system; see e.g.~\cite{Galias2021}.\llabel{rigNumEnd}

\llabel{nodeFocusBegin}Due to the topological approach, our results are limited to topological invariants. In particular, the results remain valid when the phase space is transformed by a homeomorphism, or a topologically conjugate dynamical system is considered. For example, when we find an asymptotically stable fixed point (a sink), we do not distinguish between a stable node and a stable focus.\llabel{nodeFocusEnd}

\llabel{coarseResBegin}Due to the finite resolution of our computations in the phase space and in the parameter space, this approach cannot be used directly to analyze fine structures in the phase space (like chaotic sets) or in the parameter space (like Arnold tongues or shrimp-shaped periodic regions; see \cite{De_Souza2012-vs}). However, this method effectively classifies the dynamics on the coarse level in a fully automatic manner to create a catalog of the different types of dynamics found. A possible application of this method could be to narrow down the search for bounded invariant sets that can be further investigated using other methods. For example, one could actually prove the existence of periodic orbits of specific stability type~\cite{Galias2001} or apply classical tools to track down attractors and to obtain a more precise bifurcation diagram.\llabel{coarseResEnd}

This method is robust in the sense that the results are also valid for all nearby systems. A drawback is that we only see features that are structurally stable, so we will miss very subtle phenomena. We discuss this further in Sec.~\ref{sec:final}.

We explain some details of the method in the following subsections. \llabel{high-dim-begin}Since the method and software that we use also apply without changes to the phase space and the parameter space of dimensions above $2$, we describe the method in a dimension-independent way. Indeed, this approach was successfully applied to a $4$-dimensional system in \cite[Sec.~5]{knipl-2015} or to a $6$-dimensional system in \cite[Sec.~6]{Day2004}, and analysis with $3$ parameters varying at the same time was conducted in~\cite{arai-2009}. However, visualization in higher dimensions is more demanding, so we apply it to the $2$-dimensional model \eqref{eq:model} with $2$ varying parameters in the current work.

Using an $n$-dimensional pixelation approach to the analysis of dynamics is an effective method for the discretization of the phase space that was used, for example, by M.~Dellnitz and his research group \cite{GAIO,Dellnitz1997} and by A.~Szymczak \cite{Szymczak1997}.\llabel{high-dim-end}

We would like to note that the method we apply also has an extension to dynamical systems with continuous time (flows) \cite[Sec.~5]{knipl-2015}, where it meets other challenges, e.g., related to rigorous integration of ODEs~\cite{miyaji-2016}.

\subsection{Isolated invariant sets}
\label{sec:isol}

Consider a dynamical system with discrete time generated by a continuous map $f\colon X\to X$. In case of a non-invertible map, this is formally called a semi-dynamical\llabel{semi-dynam} system; see e.g. \cite{Mrozek1989}, but we will still call it a dynamical system to keep the terminology simple.

\llabel{terminologyBegin}To fix some terminology, a point $x \in X$ such that $f^n(x)=x$ for some $n > 0$ is called a \emph{periodic point}. The smallest number $n > 0$ such that $f^n(x)=x$ is called its \emph{period}. A periodic point of period $1$ is called a \emph{fixed point}.

An \emph{orbit} of a point $x \in X$ is the set $\{f^n(x) : n \in \mathbb{Z}\}$ if $f$ is invertible or $\{f^n(x) : n \geq 0\}$ otherwise. If $x$ is a periodic point then its orbit is called a \emph{periodic orbit}.\llabel{terminologyEnd}

\begin{defin}
\label{def:inv}
We say that $S \subset X$ is an \emph{invariant set} with respect to $f$ if $f (S) = S$.
Moreover, the \emph{invariant part} of a set $N \subset X$ with respect to $f$, denoted $\inv(N,f)$, is the largest (in terms of inclusion) subset of $N$ that is invariant.
\end{defin}

Note that we do not require that an invariant set $S$ is a \emph{minimal} set such that $f(S)=S$. In particular, the union of a family of invariant sets is also an invariant set. In this way, we can talk about Morse decompositions of invariant sets discussed in the following section; see also~\cite{mischaikow2002conley}.

\begin{defin}
\label{def:isol}
A compact set $N \subset X$ is called an \emph{isolating neighborhood} with respect to $f$ if
$\inv(N, f) \subset \mathop{\textrm{int}} N$.
A set $S \subset X$ is called an \emph{isolated invariant set} if there exists an isolating neighborhood $N$ for which $\inv(N,f)=S$.
\end{defin}

We are interested in capturing all the bounded \llabel{chainRec}isolated invariant sets of the system such as fixed points, periodic orbits or chaotic sets, whether stable or unstable.
We are also interested in determining possible connections between them in the form of a Morse decomposition with the Conley index used to qualify the stability type of the Morse sets, as explained in the following subsections.


\subsection{Morse decomposition}
\label{sec:Morse}

Let us assume that there exists a sufficiently large compact set that contains all the bounded invariant sets of the system that we consider. This is the case, for example, if there exists a bounded absorbing set (see Sec.~\ref{sec:rectangleB}) and the phase space $X$ is a subset of a finite-dimensional Euclidean space. Then a sufficiently large $B\subset X$ is an isolating neighborhood such that the set $S=\inv(B,f)$ contains all the bounded invariant sets of the system, both attracting and unstable. However, knowing $B$ or even $S$ itself yields very coarse information. To paint a clearer picture of the system, $S$ can be split into a finer structure of isolated invariant sets between which the dynamics is gradient-like, as defined below.

\begin{defin}
\label{def:Morse}
Let $S$ be an isolated invariant set for $f$.
Any finite collection of its pairwise disjoint isolated invariant subsets $\{M_p\}_{p\in\mathcal{P}}$ is called a \emph{Morse decomposition} of $S$ if there is a strict partial order $\succ$ on $\mathcal{P}$ such that for all
$x\in S\setminus\bigcup_{p\in\mathcal{P}}M_p$
there exist $p, q \in \mathcal{P}$ such that $p\succ q$ and $\omega(x)\subset M_q$ and $\alpha(x)\subset M_p$, where $\omega(x)$ and $\alpha(x)$ are the corresponding limit sets of $x$; in this case, we say that \emph{there exists a connection from $M_p$ to $M_q$.}
Each of the sets $M_p$ is called a \emph{Morse set}.
\end{defin}

A Morse decomposition can be represented by means of a directed graph with vertices corresponding to the Morse sets and directed edges indicating the relation $\succ$. For clarity of presentation, we plot a transitive reduction of this graph; then $p \succ q$ if and only if there exists a directed path from the vertex representing $M_p$ to the vertex representing $M_q$ in this graph.

As a warning, we point out that the paths in the graph representing a Morse decomposition only provide partial information. If there exists a connection from $M_p$ to $M_q$, then $p \succ q$ and there exists a corresponding path. However, the opposite implication is not true. There might exist a path in the graph from the vertex corresponding to $M_p$ to the vertex corresponding to $M_q$ even if there exists no connection from $M_p$ to $M_q$. On the other hand, if there is no directed path between certain vertices in the graph, then it is true that there exists no connection between the corresponding Morse sets.

\llabel{aimSetsBegin}In this approach, we do not aim to find the precise location of invariant sets, but rather to estimate their location in terms of isolating neighborhoods constructed at a prescribed resolution. In contrast to simple numerical simulations, we find all stable and unstable sets. Let us briefly introduce this method in our setting.\llabel{aimSetsEnd}


\subsection{Numerical Morse decomposition}
\label{sec:numMorse}

Assume $X=\mathbb{R}^n$, and a cuboid $B = \Pi_{i=1,\ldots,n} [a_i,b_i] \subset \mathbb{R}^n$ is an isolating neighborhood with respect to a continuous map $f \colon X \to X$. Consider a uniform rectangular grid $\mathcal{G}(B)$ in $B$.

Representing the map $f$ restricted to $B$ by a directed graph with vertices $\mathcal{G}(B)$ and edges reflecting the map $f$ makes it possible to find a \emph{numerical Morse decomposition} consisting of \emph{numerical Morse sets}, built as finite unions of grid elements, that are isolating neighborhoods of Morse sets in a real Morse decomposition of $\inv(B,f)$, as in Def.~\ref{def:Morse}. An effective method for this purpose, based on fast graph algorithms, was introduced in \cite{arai-2009} and later developed further and expanded; see, e.g., \cite{pilarczyk-2023}. The involvement of computational topology tools for the computation of homomorphisms induced in homology
and the Conley index in this context \cite{pilarczyk-2008} provides additional information on the stability of the constructed numerical Morse sets.

\llabel{numMorse}We remark that we do not actually construct Morse sets, but rather their isolating neighborhoods; therefore, we call them \emph{numerical} Morse sets.


\subsection{Conley index}
\label{sec:Conley}

Given a numerical Morse set $M$, the graph representation of the map $f$ is used to compute a rigorous outer bound on $f(M)$ in terms of a union of grid elements; let us denote it by $F(M)$. Then the pair of compact sets $P_1 := M \cup F(M)$ and $P_2 := \cl (F(M) \setminus M)$ satisfies the conditions for an index pair for the Conley index that provides qualitative information on the stability of $\inv(M,f)$. The set $P_2$ is called the \emph{exit set} of $M$ because all trajectories can only exit the isolating neighborhood $M$ through $P_2$. The relative homology of $(P_1,P_2)$ together with the homomorphism induced by $f$ in homology are the core ingredients of the \emph{Conley index} of $M$.

Without delving into technical details, let us point out two key features of the Conley index and explain what kind of information this index provides. First, the Conley index does not depend on the choice of an index pair as long as $S := \inv(\cl(P_1\setminus P_2),f)$ remains the same; therefore, the Conley index is a property of an isolated invariant set, not of a specific index pair or isolating neighborhood. Second, if the Conley index is nontrivial, then $S \neq \emptyset$; this fact can be used to 
prove that we have captured a nontrivial invariant set. Since all computations are performed using rigorous numerical methods, this is actually a computer-assisted proof of this fact. We refer to \cite{mischaikow2002conley} for a brief introduction to Conley index theory and to \cite{KMM2004} for index pairs and the Conley index in the context of cubical sets.

At this point, we would like to draw attention to the fact that the Conley index is a topological invariant that reflects the features of the dynamics inside $M$ as they look from the outside; therefore, the actual dynamics inside $M$ can be different or considerably more complex. The Conley index of a hyperbolic periodic point implies the existence of a periodic point in $S$, which this follows, e.g., from \cite[Thm.~4]{Mrozek1989}\llabel{hypPointExist} or \cite[Thm.~3]{Mrozek1990} through the fixed-point index theory. However, the actual stability of that periodic point might be different from what could be inferred from the Conley index, and the set $S$ might have a more complex structure than just a hyperbolic periodic point.

\begin{figure}[htbp]
\centering
\includegraphics[scale=0.8]{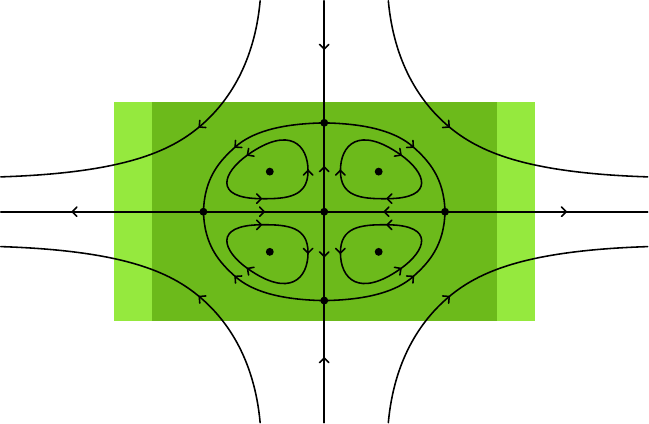}
\caption{Example of an index pair whose Conley index is the same as the Conley index of a hyperbolic fixed point, but the corresponding isolated invariant set (see Def.~\ref{def:inv}) actually consists of several smaller invariant subsets. The isolating neighborhood is shown in green, together with the corresponding exit set drawn in a lighter tone. The trajectories are shown for a flow; consider the time-$1$ map for a discrete-time dynamical system.}
\label{fig:conleyExample}
\end{figure}

\llabel{fig2begin}Figure~\ref{fig:conleyExample} gives an example of the purpose and limitations of this method. The isolating neighborhood constructed at the given resolution has the Conley index of an unstable fixed point (a saddle). This provides rigorous (mathematically proven) bounds for an isolated invariant set and the information that the set is non-empty and contains a fixed point. However, we do not obtain any more specific information about the actual dynamics on a finer scale. In particular, this method does not provide the characterization of individual solutions.\llabel{fig2end}

To sum up, the most important conclusion that we draw from a nontrivial Conley index is the existence of a nontrivial isolated invariant set inside the isolating neighborhood found. Indeed, it can be argued that from the point of view of applications, the actual fine-scale microscopic structure of the isolated invariant set found is of lesser importance than the fact of its non-emptiness and the effective bound for its location provided by the isolating neighborhood; we refer to \cite{luzzatto-2011} for a more in-depth discussion of this question.

\subsection{Conley--Morse graph}
\label{sec:CMGraph}

A graph that represents a Morse decomposition together with topological information provided by the Conley index of each numerical Morse set is called a \emph{Conley--Morse graph}, or a \emph{CM graph} for short.

To explain the intuition behind the information provided by the Conley index in a somewhat informal way, one can say that this topological information consists of two ingredients: \llabel{conIndBegin}rough characterization of the topology of the isolating neighborhood in terms of its homology groups\llabel{conIndEnd} and information on its stability. The first ingredient is obtained at the homological level and concerns the static information, so it only reflects the number of connected components and the number of holes of different dimensions, such as the hole in the circle, the void in the sphere, etc. The stability information essentially reduces to indicating the number of unstable directions relative to the dimension of the embedding space and whether the map $f$ preserves or reverses the orientation in the unstable directions considered collectively.

In our case of a planar system, we typically have an attractor (a sink; $0$ unstable and $2$ stable directions), a hyperbolic saddle ($1$ unstable and $1$ stable direction) or a repeller (a source; $2$ unstable and $0$ stable directions), with a map that either preserves or reverses the orientation. We have designed pictograms to schematically indicate this information. We explain and illustrate them in Fig.~\ref{fig:symbols} and refer to \cite[Figs. 2--4]{pilarczyk-2023} for phase space portrait examples corresponding to many of these cases.

\llabel{trivialIndex}We emphasize the fact that our method first finds isolating neighborhoods and then computes their Conley indices. Therefore, it may happen that an isolating neighborhood is constructed where there is a slow-down in the dynamics only, for example, near a bifurcation, as shown in Fig.~\ref{fig:sn_bifur}, and the invariant part of this neighborhood is in fact empty. Obviously, in this case the Conley index is trivial. The last pictogram in Fig.~\ref{fig:symbols} shows such a situation.

\llabel{pict_expr}We are aware of the fact that the power of expression of these pictograms is limited, and additional pictograms may be needed for other types of invariant sets found (such as invariant tori, for example), but we feel that they provide the most relevant information in a clearer way than expressing the homology of the index pair and the index map using formulas, as was done so far (see, e.g., \cite{arai-2009,pilarczyk-2023}).

\begin{figure}[htbp]
\centering
\small
\begin{tabular}{cp{7cm}}
\raisebox{-0.5\totalheight}{\includegraphics[scale=1]{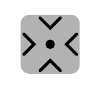}} & a stable (attracting) fixed point \\
\raisebox{-0.6\totalheight}{\includegraphics[scale=1]{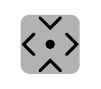}} & a saddle fixed point with an orientation-preserving map \\
\raisebox{-0.8\totalheight}{\includegraphics[scale=1]{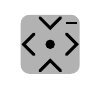}} & a saddle fixed point with a flip in the unstable direction indicated by the ``$-$'' symbol in the corner (this is a saddle fixed point inscribed in a M\"{o}bius band if one thinks to a suspension of the map) \\
\raisebox{-0.5\totalheight}{\includegraphics[scale=1]{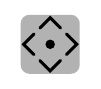}} & an unstable (repelling) fixed point \\
\raisebox{-0.5\totalheight}{\includegraphics[scale=1]{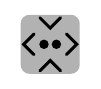}} & the orbit of a saddle period-$2$ point \\
\raisebox{-0.6\totalheight}{\includegraphics[scale=1]{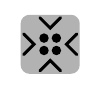}} & the orbit of a stable (attracting) periodic point of period $4$ \\
\raisebox{-0.6\totalheight}{\includegraphics[scale=1]{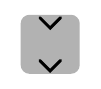}} & the trivial (empty) set with trajectories apparently going through its isolating neighborhood \\
\end{tabular}

\caption{We represent different types of the Conley index by pictograms. The arrows indicate the stability type (source, sink, saddle). The ``$-$'' symbol in the corner indicates a flip in the unstable direction. The dots inside indicate the number of connected components of the isolating neighborhood.}
\label{fig:symbols}
\end{figure}

\subsection{Continuation}
\label{sec:continuation}

A graph representation of the map $f$ and thus the numerical Morse decomposition is computed with rigorous numerical methods based on \emph{interval arithmetic} \cite{Moore2009}. The idea behind this technique is to perform computations on intervals instead of numbers in such a way that all rounding errors are taken into account (by rounding the results downwards or upwards at the endpoints of intervals) and the resulting intervals \emph{contain} exact results. In other words, we compute an outer estimate for the resulting intervals each time. This inevitably leads to overestimates, so we work with \emph{inclusions} instead of \emph{equations}, one could say, and therefore we need to use topological tools to draw mathematically valid nontrivial conclusions. One can compare this kind of reasoning to using Bolzano's theorem to conclude the existence of a root of a continuous function from inequalities on the endpoints of the interval.

Assume $\Lambda$ is a product of compact intervals of parameters. Consider a map that depends on these parameters, that is, $f \colon X \times \Lambda \ni (x, \lambda) \mapsto f_{\lambda}(x) \in X$. Introduce a rectangular grid $\mathcal{G}(\Lambda)$ in $\Lambda$. By replacing individual parameters in the formula for $f$ with intervals of parameters, we can carry out the procedure of constructing the CM graph described above for each $L \in \mathcal{G}(\Lambda)$. Since interval arithmetic is used throughout this process, the CM graph constructed for $L$ is valid for $f_{\lambda}$ for every $\lambda \in L$, \llabel{inclusions}specifically because the conditions for isolating neighborhoods and index pairs are formulated in terms of inclusions of sets that can hold true even if some overestimates are made in the computations.

If the numerical Morse sets in CM graphs constructed for two cells of parameters $L_1, L_2 \in \mathcal(G)(\Lambda)$ such that $L_1 \cap L_2 \neq \emptyset$ are in one-to-one correspondence, then the dynamics found for parameter values in $L_1$ and in $L_2$ is equivalent. The equivalence classes of this relation are called \emph{continuation classes}.

We consider a list of continuation classes with the corresponding CM graphs to be a complete description of global dynamics found in the system at the finite resolution imposed by the grid $\mathcal{G}(X)$ across the parameter set $\Lambda$. If $\Lambda$ is a product of two intervals, then the computed continuation classes can be easily visualized, as we show in Sec.~\ref{sec:results}; see Fig.~\ref{fig:contDiag}.

\section{Results and discussion}
\label{sec:results}

Before applying the CMAD method introduced in Sec.~\ref{sec:method} to the Andrecut–Kauffman model \eqref{eq:model}, we have to fix bounded ranges of parameters that form $\Lambda$ and a rectangular grid $\mathcal{G}(\Lambda)$ (see Sec.~\ref{sec:continuation}), as well as an isolating neighborhood $B \subset \mathbb{R}^2$ in which to perform the computations, and a rectangular grid $\mathcal{G}(B)$ (see Sec.~\ref{sec:numMorse}) to fix a finite resolution at which we perceive the system. We provide and justify our choices in Sec.~\ref{sec:param} and Sec.~\ref{sec:rectangleB}. In the next sections, we describe and discuss the results obtained: We show different CM graphs that we found in the main continuation classes and explain the information on the dynamics that they provide, including possible bifurcations inferred from the adjacency of continuation classes.
Additionally, in Sec.~\ref{sec:resol}, we discuss the impact of the phase space resolution on the accuracy of the dynamics that can be observed.

\subsection{\texorpdfstring{The parameter set $\Lambda$}{The parameter set}}
\label{sec:param}

The topological--numerical method and the corresponding software are capable of performing computations for a multi-dimensional space of parameters. Therefore, in principle, we could take wide ranges of $\alpha_1$, $\alpha_2$, $\beta_1$, $\beta_2$, $\varepsilon$, and perform computations for several values of $n$. However, the cost of such computations would be prohibitive and the $5$-dimensional continuation classes would be difficult to analyze and visualize. Therefore, for the sake of simplicity, we chose to fix $\beta = \beta_1 = \beta_2 = 0.2$, $\varepsilon = 0.8$, and $n = 3$ and consider the two-dimensional space of varying parameters $(\alpha_1, \alpha_2) \in \Lambda = [0,80] \times [0,80]$. Based on the literature and our numerical experiments, these seem to be the most relevant parameters that additionally make it possible to observe a variety of different dynamical behaviors. The chosen ranges of $\alpha_1$ and $\alpha_2$ cover most of the values of these parameters considered in the literature (see Sec.~\ref{sec:state}). The selected level of $\beta$ appears in other research focusing on this system, and $\varepsilon = 0.8$ was selected due to the noticeable bistability of the system in the range $(0.7,1)$ found in \cite{rosman-2025}. \llabel{n3Begin}Our choice of $n=3$ is justified by the fact that the dynamics for $n=2$ is considerably simpler, and $n=3$ was also chosen in some literature (e.g.~\cite{De_Souza2012-vs}).\llabel{n3End} We subdivide the chosen parameter space $\Lambda$ into the uniform rectangular grid $\mathcal{G}(\Lambda)$ of $160 \times 160$ squares of equal size.

\subsection{\texorpdfstring{The bounded isolating neighborhood $B$}{The bounded isolating neighborhood B}}
\label{sec:rectangleB}

We chose to conduct the computations in the rectangle $(x,y) \in B = [0,101] \times [0,101]$.
The set $B$ was divided into a rectangular grid of $4096 \times 4096$ rectangles.
This subdivision defines the actual resolution at which the graph representation of the map was computed and the numerical Morse decomposition was constructed.
In fact, the software uses a gradual refinement technique to efficiently narrow down its search for $\inv(B,f)$ within $B$ (see \cite[Sec.~4.2]{arai-2009}), so it starts the computations at a coarser resolution and then halves the grid in each direction a few times to continue the computations at the finer resolutions. This is the reason why the number of subdivisions of $B$ in each direction is a power of $2$.

This selection of the set $B$ is justified by the result proved in \cite[Proposition~4.2]{rosman-2025}, using the concepts introduced in \cite{pilarczyk-graff-2024}, that all the bounded invariant sets of the Andrecut--Kauffman model are contained within the set $[0,b_1]\times[0,b_2]\subset\mathbb{R}^2$, where $b_i := (\alpha_i + c_i)/(1-\beta)$ with $i\in\{1,2\}$ and arbitrarily small $c_i>0$. For the chosen values of $\alpha_1$, $\alpha_2$, and $\beta$, the absorbing set with $c_1$ and $c_2$ approaching $0$ would be $[0,100]\times[0,100]$. A slightly larger $B$ was used to account for small nonzero values of $c_1$ and $c_2$.

\subsection{Numerical computations}
\label{sec:comput}

For the purpose of computing CM graphs and continuation classes for the dynamical system induced by \eqref{eq:model} with the fixed parameters and $\Lambda$ described in Sec.~\ref{sec:param} and the phase space restricted to $B$ defined in Sec.~\ref{sec:rectangleB}, we used a software implementation of the method described in Sec.~\ref{sec:method} available on \url{https://www.pawelpilarczyk.com/gen2d/}. This software uses the parallelization scheme introduced in \cite[Sec.~3]{pilarczyk-2010} and runs efficiently on a computer cluster or on a computer with a CPU that features multiple cores.

The computations were performed on a cluster of computers with Intel\textregistered{} Xeon\textregistered{} Platinum 8268 processors clocked at 2.90GHz. The total time used was almost 15 CPU hours. The results of these computations and an interactive browser of CM graphs and numerical Morse sets are available on \url{https://www.pawelpilarczyk.com/gen2d/}.

\llabel{attrSim}In addition to applying the CMAD method that constructs isolating neighborhoods at a fixed resolution, we performed classical numerical simulations in selected cases to obtain a more precise approximation of the attractors observed. Since the isolating neighborhoods provide validated results on the location of all types of invariant sets, in the case of attractors numerical simulations yield considerably finer approximations of them (although this is not a computer-assisted proof). The computed attractors are overlaid in black on top of the computed isolating neighborhoods in all the figures depicting numerical Morse sets in the phase space.
These attractor simulations were performed for a uniform grid of $512 \times 512$ initial conditions in $B$. For each of these points, only its position after $10\,000$ iterations has been plotted. The simulations were performed in parallel for all initial conditions using NVIDIA\textregistered{} CUDA\texttrademark{} technology on a single NVIDIA\textregistered{} RTX3060 6GB GPU and took approximately $1.1$ seconds per figure.

\subsection{Continuation diagram}
\label{sec:contDiag}

The resulting continuation diagram is shown in Fig.~\ref{fig:contDiag}. Let us explain how this diagram is created and what information it provides.

\begin{figure}[tbp]
    \centering
    \includegraphics[width=\linewidth]{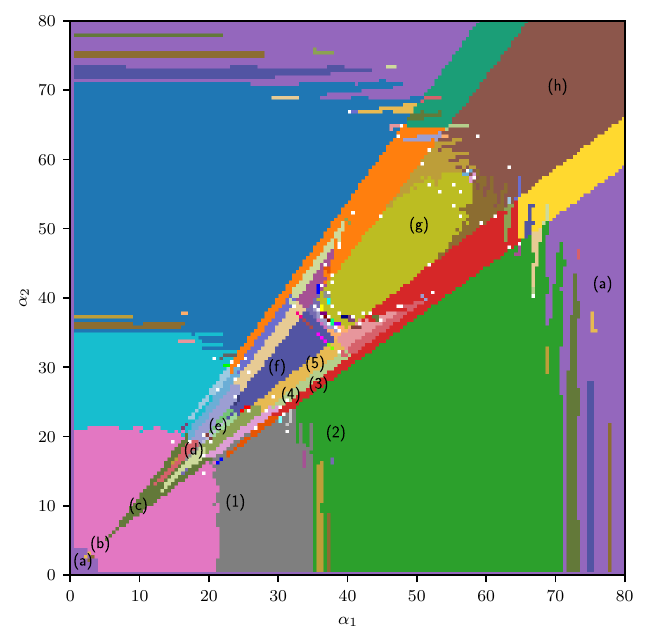}
    \caption{Continuation diagram described in Sec.~\ref{sec:contDiag} for the Andrecut--Kauffman model \eqref{eq:model} with $(\alpha_1,\alpha_2) \in [0,80] \times [0,80]$, $\beta=0.2$, $\varepsilon=0.8$, and $n=3$. Each continuation class (see Sec.~\ref{sec:continuation}) is shown in a different color (except for classes containing exactly one element, all shown in white). Labels (a)--(h) and (1)--(5) indicate classes discussed in Sec.~\ref{sec:diag} and Sec.~\ref{sec:different}, respectively. The colors of continuation classes shown in this diagram are also used in Fig.~\ref{fig:lyapunovDiag}.}
    \label{fig:contDiag}
\end{figure}

\llabel{diagExplBegin}The set of parameters $(\alpha_1, \alpha_2) \in [0,80] \times [0,80]$ is divided into $160 \times 160$ squares of equal size.
Each square corresponds to a subset of parameters.
For example, the square in the bottom left corner corresponds to the parameters $(\alpha_1, \alpha_2) \in [0,0.5] \times [0,0.5]$.
The numerical Morse decomposition constructed at finite resolution in the phase space for each square of parameters and the corresponding CM graph are rough\llabel{rough} enough to be valid for all the parameters in the entire square of parameters.

CM graphs constructed for pairs of adjacent parameter squares are compared with each other to see if they represent the same information about the dynamics, as explained in Sec.~\ref{sec:continuation}.
If the comparison is affirmative, we say that the adjacent squares belong to the same continuation class.
We extend this relation by transitivity to split the entire set of $160 \times 160$ squares into disjoint continuation classes.
We then choose a different color for each continuation class consisting of more than one element and we plot these squares using the chosen colors while leaving the singletons white.

As a result, the set of parameters is painted in such a way that the perception of system's dynamics is qualitatively the same for all parameters in adjacent squares of the same color in the diagram.\llabel{diagExplEnd}

There are $110$ continuation classes, about half of which are built of more than $10$ parameter squares each.
Note the symmetry in the diagram with respect to the line $\alpha_1 = \alpha_2$.
This is due to the symmetry in the equations when taking $\beta_1 = \beta_2$ in our case.
The phase portraits on the opposite sides of this line are symmetric reflections of each other with respect to the line $x = y$.
In what follows, we restrict ourselves to describing the dynamics for $\alpha_1 \geq \alpha_2$.

In the continuation diagram, two groups of regions can be distinguished: regions located close to the diagonal of the diagram and regions significantly distant from the diagonal.
The regions located near the diagonal are not large but more numerous.
These regions become progressively larger as the values of the parameters $\alpha_1$ and $\alpha_2$ increase.
In contrast, the regions far from the diagonal are considerably larger, but there are fewer of them.
Therefore, we will first discuss the dynamics found in selected regions along the diagonal, labeled (a)--(h) in Fig.~\ref{fig:contDiag}, corresponding to increasing the parameters $\alpha_1$ and $\alpha_2$ simultaneously; see Sec.~\ref{sec:diag}.
Next, we will examine what kinds of dynamics can be seen when the parameter values are different but gradually approach the diagonal, which corresponds to regions (1)--(5) in Fig.~\ref{fig:contDiag}; see Sec.~\ref{sec:different}.

But before that, we will first show and briefly discuss the largest Lyapunov exponent and the sequence of bifurcations observed along the diagonal ($\alpha_1 = \alpha_2$) in numerical simulations (Sec.~\ref{sec:Lyapunov}), and say a word of caution regarding bifurcations observed using the CMAD method (Sec.~\ref{sec:bifurcations}).

\subsection{Bifurcation diagram along the diagonal and the largest Lyapunov exponent}
\label{sec:Lyapunov}

Most previous work on the Andrecut--Kauffman model \eqref{eq:model} focused primarily on using numerical simulations to compute bifurcation diagrams and the largest Lyapunov exponent, and partitioning the parameter space into chaotic and periodic regions \cite{andrecut_main,rosman-2025,subramani-2023}. Those methods provide valuable information about the types of behavior expected for different parameter values. However, in the context of physical systems where noise is common and measurements have limited precision, exact location and shape of attractors obtained from simulations can be less relevant than general knowledge about shape and location of isolating neighborhoods. The CMAD method provides such a coarse bird's-eye view of the dynamics present in the system. It categorizes parameter values into classes based on numerical Morse decompositions consisting of isolating neighborhoods found and their Conley indices, highlighting global behavior while mostly ignoring small-scale phenomena. Additionally, it captures unstable invariant sets, such as repellers (sources) and saddle points, which are largely invisible to standard numerical simulations.

It is natural to ask whether there is a relation between the results that can be obtained in numerical simulations (such as the largest Lyapunov exponent and the bifurcation diagram) and the continuation classes that we found. For this purpose, we computed the largest Lyapunov exponent for several values of the parameters along the diagonal of the continuation diagram and an attractor of the system for the same parameters.
The same set of 128 randomly generated initial conditions sampled from $\left\{(x, y)\in[0.01, 50]\times[0.01, 50]\colon x\geq y\right\}$ using the uniform distribution were used for each parameter value in the numerical simulations. The plotted values were estimated from 1000 iterations of the system, following 5000 iterations of burn-in.

\begin{figure}[htbp]
    \centering
    \includegraphics[height=5cm]{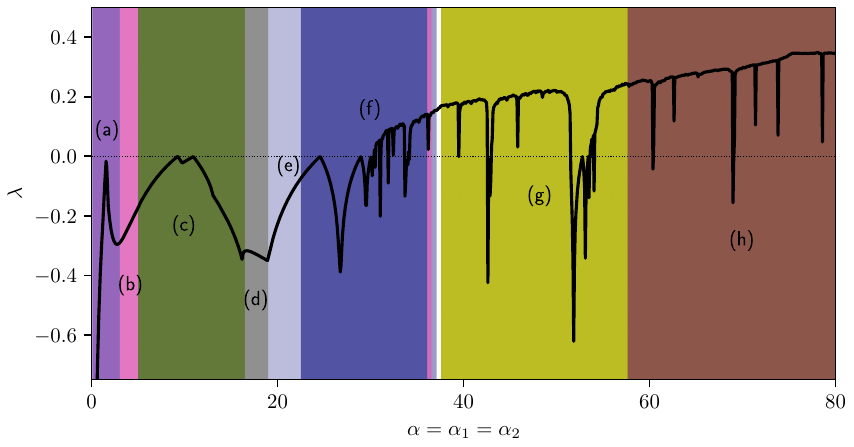}\hspace*{12pt} \\
    \hspace*{7pt}\includegraphics[height=5cm]{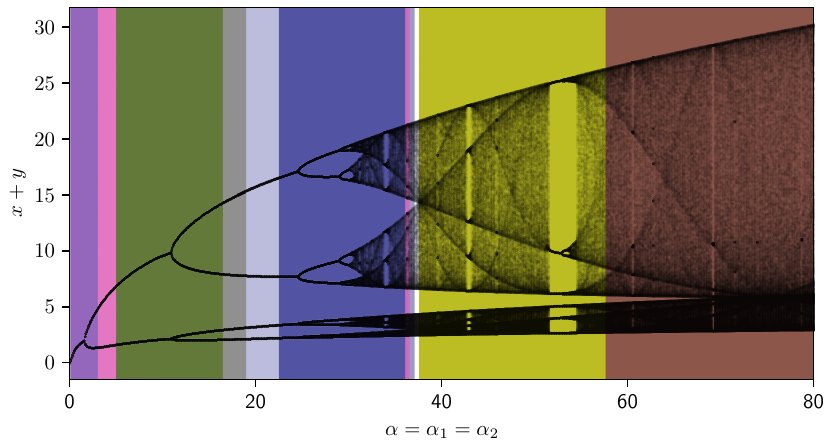}\hspace*{12pt}
    \caption{The largest Lyapunov exponent (top) and the bifurcation diagram (bottom) calculated for $\alpha_1=\alpha_2$ for 1024 equally spaced parameters with initial conditions in the phase space restricted to $x \geq y$ due to the symmetry in the system. Parameter values and colors correspond to those in Fig.~\ref{fig:contDiag}.}
    \label{fig:lyapunovDiag}
\end{figure}

The results of our numerical simulations are gathered in Fig.~\ref{fig:lyapunovDiag}. The largest Lyapunov exponent and the bifurcation diagram (showing the sum of both coordinates in the phase space on the attractor) are plotted as functions of $\alpha = \alpha_1 = \alpha_2$.
The colors used in the background of the plots match the continuation classes shown in Fig.~\ref{fig:contDiag} through which the diagonal $\alpha_1 = \alpha_2$ traverses.

\llabel{contBifBegin}The computed continuation classes represent equivalent dynamics found on the level of isolating neighborhoods. Moving to another class corresponds to a change in the CM graph, which is often related to an actual bifurcation in the system, as discussed in Sec.~\ref{sec:bifurcations}.

An important feature of the CMAD method is that it is not limited to stable sets, but it also captures bifurcations of unstable sets; such bifurcations are not shown in the bifurcation diagram obtained in numerical simulations. For example, the first period doubling corresponding to the transition from the violet region on the very left-hand side of the graphs to the pink region is followed by a pitchfork bifurcation of the saddle fixed point in the perpendicular direction when entering the green region; we discuss it further in Sec.~\ref{sec:regionC}.

Some other transitions between classes may reflect more subtle changes or may correspond to the level of detail of dynamics at the fixed resolution. For example, the transition from the dark blue region (containing $\alpha = 30$) to the light green region (containing $\alpha = 40$) goes through three intermediate regions, beginning with a change in the connections between the isolating neighborhoods in the CM graph; we refer to the online browser of the results for details.

Note that period doubling cascades and chaotic dynamics with periodic windows are often located within common continuation classes, as shown in Fig.~\ref{fig:lyapunovDiag} for $\alpha \geq 22.5$. This is because large isolating neighborhoods are created for such attractors that are continued through to adjacent parameter boxes. Therefore, as one should expect, the results obtained at the coarse limited resolution reflect large-scale changes in the dynamics like some period doubling bifurcations or pitchfork bifurcations of periodic points, but cannot be directly used to detect the presence of chaos in the system; the latter can be determined, for example, by the existence of a positive Lyapunov exponent. Indeed, the presence of chaotic dynamics causes the CMAD method to construct a large isolating neighborhood that contains the entire chaotic attractor, without revealing the details of the dynamics inside this neighborhood. In this way, numerical simulations are useful to complement the topological analysis.\llabel{contBifEnd}

\subsection{Bifurcations}
\label{sec:bifurcations}

When using the CMAD method, the analysis of dynamics and bifurcations is based only on the constructed isolating neighborhoods and their Conley indices, which provide rough information about the dynamics at a finite resolution. As shown in Fig.~\ref{fig:conleyExample} discussed in Sec.~\ref{sec:Conley}, the actual dynamics inside each isolating neighborhood may be more complex \llabel{moreComplex}than it appears when looking at the shape of the isolating neighborhood only and its stability.

In the finite resolution setting, bifurcations are perceived in a specific way. Although it is often easy to associate changes between CM graphs to the underlying bifurcations, like period doubling, the actual bifurcations in the system typically occur at nearby values of parameters, not exactly on the boundary between adjacent continuation classes. Some bifurcations may correspond to a sequence of two or more consecutive changes in the CM graph.

\begin{figure}[htbp]
    \centering
    \includegraphics[scale=0.8]{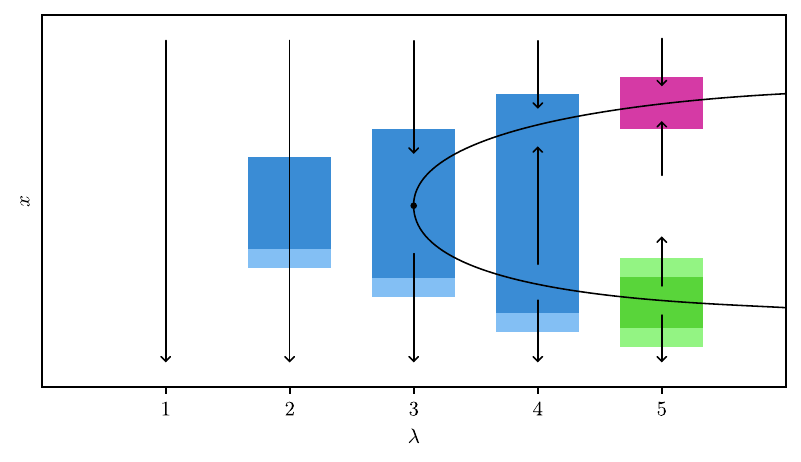}
    \caption{A saddle--node bifurcation in a 1-dimensional system seen through the lens of isolating neighborhoods at finite resolution, discussed in Sec.~\ref{sec:bifurcations}. The phase space is shown vertically, and the bifurcation takes place when the parameter $\lambda$ of the system is increased. Exit sets are shown in pale colors.}
    \label{fig:sn_bifur}
\end{figure}

Let us illustrate this feature with an example of a saddle--node bifurcation in a one-dimensional map that takes place as the parameter $\lambda$ of the map increases. Figure~\ref{fig:sn_bifur} shows isolating neighborhoods constructed for a few increasing values of the parameter $\lambda$. In this example, five qualitatively different situations occur. At $\lambda=1$, all trajectories flow along the observed section of the phase space and no isolating neighborhood is constructed. At $\lambda=2$, the system is approaching a bifurcation and the dynamics begins to slow down. This change is picked up by the finite resolution analysis, resulting in an isolating neighborhood with the trivial Conley index and the empty invariant part. This coarse perception of the dynamics does not change as the bifurcation takes place at $\lambda=3$, although the invariant part is no longer empty. For larger values of $\lambda$, there are two fixed points, one attractor and one repeller. However, for $\lambda=4$, these two fixed points are still too close to each other to be distinguished at the selected resolution. When the distance between the two points becomes larger, the full picture can be resolved, which is shown at $\lambda=5$, where two distinct isolating neighborhoods can be seen.

Independent of the problems discussed above, when sufficiently small isolating neighborhoods are constructed, the obtained results represent a significant contribution to understanding the dynamics in applications where precision is limited and noise or disturbances are expected in the real system.
However, when a large isolating neighborhood is constructed for some parameter values, it provides very limited information about the dynamics beyond a large (little informative) bound on the location of the corresponding isolated invariant set and the information that can be derived from its Conley index.

\subsection{Dynamics along the diagonal}
\label{sec:diag}

Let us describe the dynamics of the system using the continuation classes that appear in Fig.~\ref{fig:contDiag} along the diagonal $\alpha_1 = \alpha_2$.
For each continuation class, we show its CM graph and the numerical Morse sets in the phase space cropped to a neighborhood of these sets. Additionally, we plot the exit set (if non-empty) of each numerical Morse set with a brighter color corresponding to the dark color used for the set itself.

\subsubsection{A globally stable fixed point in region (a)}
\label{regionA}

\begin{figure}[htbp]
    \centering
    \begin{subfigure}[c]{0.25\linewidth}
        \centering
        \includegraphics[scale=1]{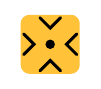}
    \end{subfigure}
    \hspace{0.05\linewidth}
    \begin{subfigure}[c]{0.5\linewidth}
        \centering
        \includegraphics[scale=0.5]{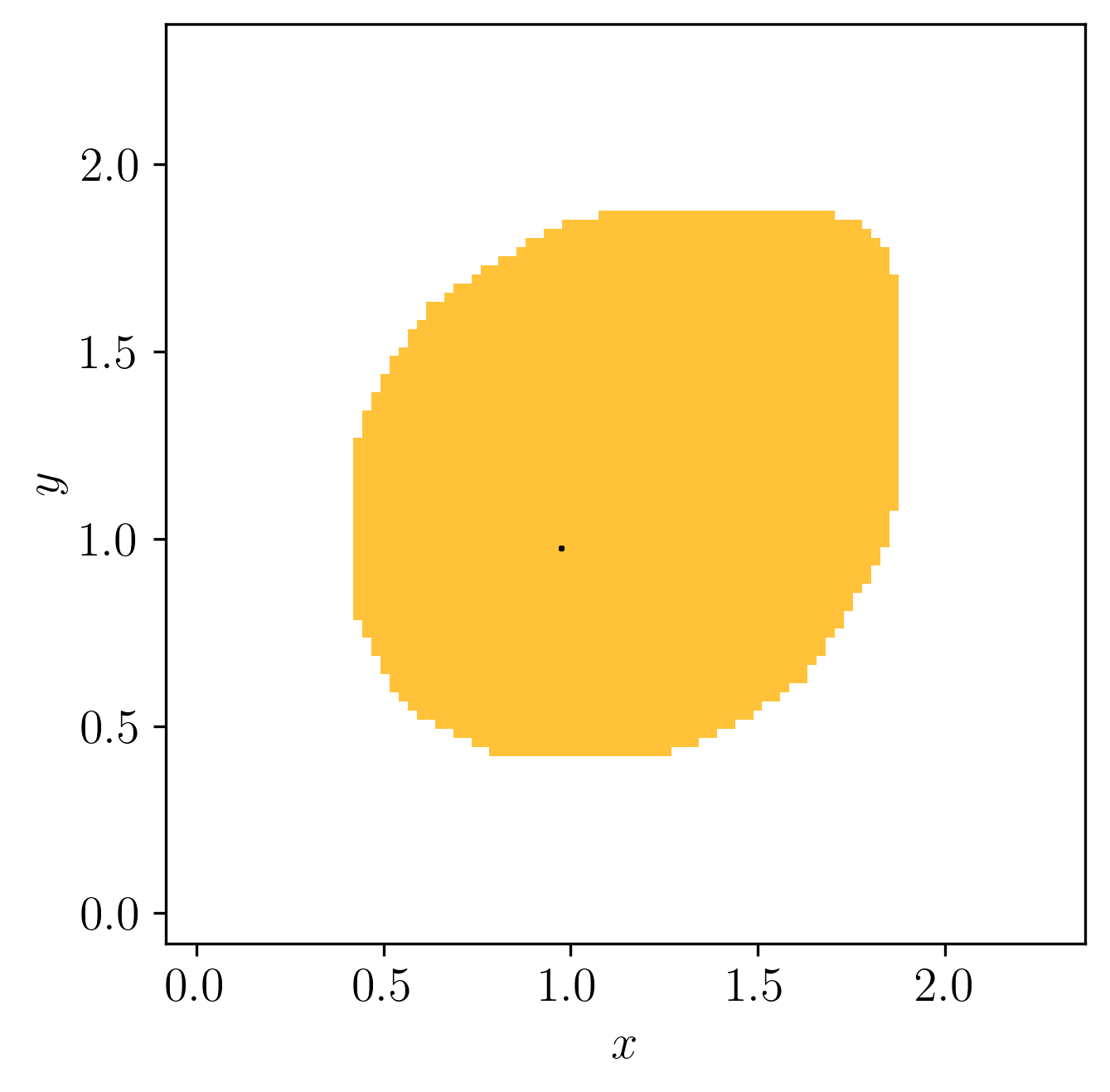}
    \end{subfigure}
    \caption{Conley--Morse graph and the isolating neighborhood found for $\alpha_1,\alpha_2\in[1.5,2]$ in region (a).}
    \label{fig:3_3}
\end{figure}

Figure~\ref{fig:3_3} shows the CM graph and the corresponding isolating neighborhood (numerical Morse set) $M$ in the phase space found for a selected parameter box taken from region (a) containing the segment $\alpha_1 = \alpha_2 \in [0,3]$.

The exit set is empty, which implies that the inclusion $f(M) \subset \interior M$ has been proven, due to the way a representation of the map was constructed (see \cite{arai-2009,pilarczyk-2023}). Since $B$ is an absorbing set for the system, this effectively proves that all trajectories enter $M$ and do not leave it. The Conley index indicates a single connected component and stability in both directions.

Based on the above information, we conclude that the dynamics in region (a) looks like a system with a global attractor that at the finite resolution at which the computations were conducted looks like a stable fixed point. The CMAD method provides its \llabel{roughFig7}rough location in the space. Numerical simulations allowed us to show the location of this fixed point more precisely; it is indicated with the black dot inside the isolating neighborhood in Fig.~\ref{fig:3_3}.

Note that region (a) includes all the parameter boxes that contain $\alpha_1 = 0$ or $\alpha_2 = 0$. For these parameter boxes, the Conley index could not be computed due to the lack of isolation in the chosen subset $B$ of the phase space (the set $M$ intersected the boundary of $B$).

\subsubsection{A stable periodic orbit of period $2$ in region (b)}
\label{sec:regionB}

\begin{figure}[htbp]
    \centering
    \begin{subfigure}[c]{0.2\linewidth}
        \centering
        \includegraphics[scale=1]{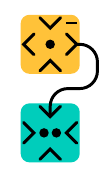}
    \end{subfigure}
    \hspace{0.05\linewidth}
    \begin{subfigure}[c]{0.5\linewidth}
        \centering
        \includegraphics[scale=0.5]{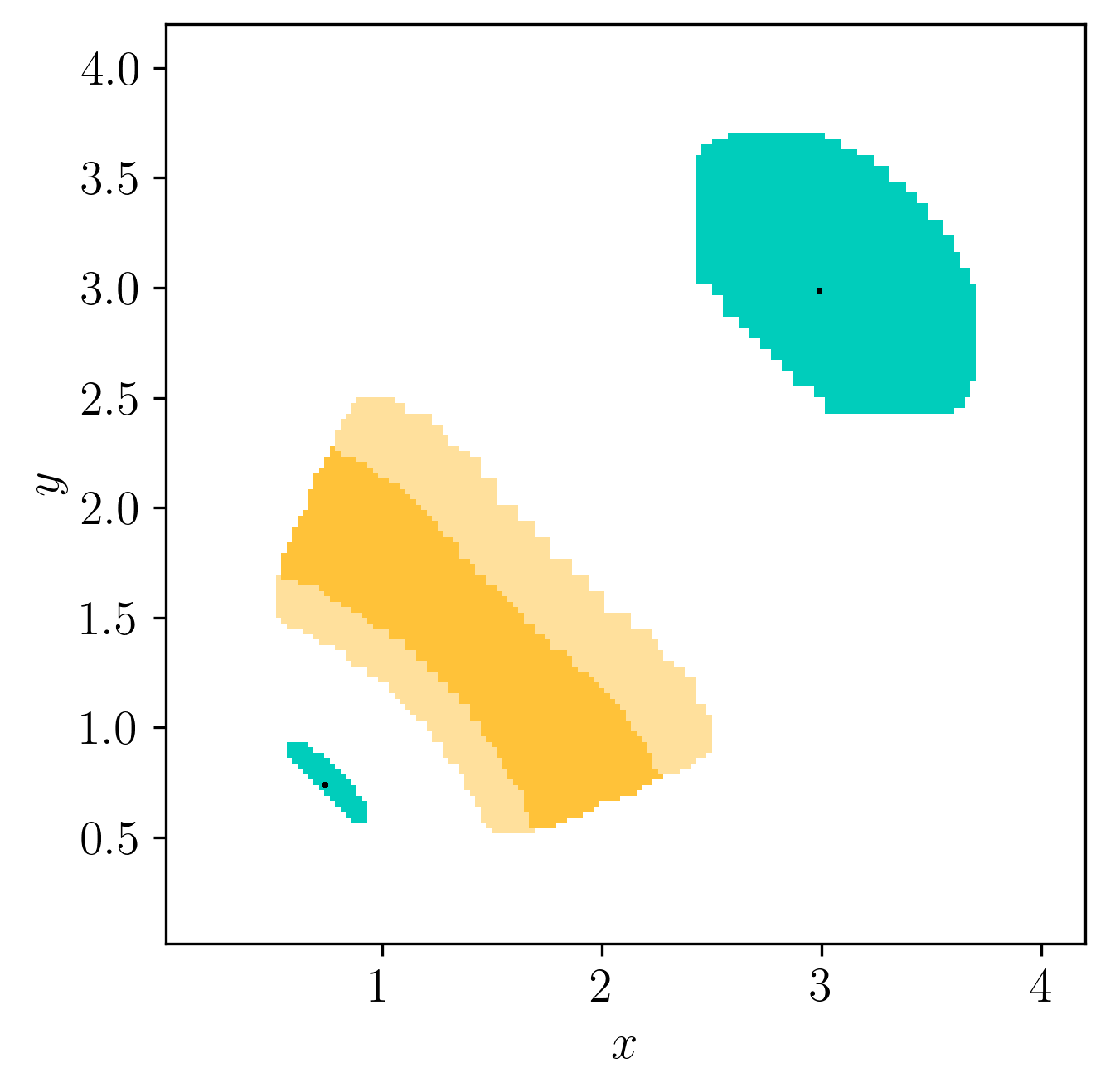}
    \end{subfigure}
    \caption{Conley--Morse graph and the isolating neighborhoods along with their exit sets (shown in pale colors) found in the phase space for $\alpha_1,\alpha_2\in[4,4.5]$ in region (b).}
    \label{fig:8_8}
\end{figure}

Increasing the values of the parameters $\alpha_1 = \alpha_2$ leads to region (b) containing the segment $\alpha_1 = \alpha_2 \in [3,5]$.
The CM graph found in that region and an example of the corresponding numerical Morse sets in the phase space are shown in Fig.~\ref{fig:8_8}.

Listed at the top of the CM graph, there is a numerical Morse set (shown in yellow) whose Conley index corresponds to a saddle fixed point.
This is also visible in the exit set of this Morse set (shown in pale yellow)\llabel{exitYellow} that does not surround this Morse set, but its two components stick out towards the other Morse set.
An edge in the CM graph connects this set to an isolating neighborhood consisting of two connected components with the Conley index of a stable periodic orbit of period~$2$.
The image of the small connected component is contained in the interior of the larger one, and \textit{vice versa}.

Transition from region (a) to region (b) corresponds to the classical period doubling bifurcation in which a stable fixed point loses its stability and a stable period-$2$ orbit appears in its vicinity. Note that just after the bifurcation, the period-$2$ orbit is located in the phase space very close to the fixed point. It can only be separated in a distinct isolating neighborhood when both invariant sets move apart well enough in the phase space, as explained in Sec.~\ref{sec:bifurcations}. A look at the bifurcation diagram (Fig.~\ref{fig:lyapunovDiag}) makes it clear that the actual bifurcation takes place around $\alpha = 1.5$, but becomes visible with our \llabel{withOurApproach}approach only after $\alpha$ has been increased beyond $\alpha = 3$. These specific values come from the discussion in Sec.~\ref{regionA}. Region (a) corresponds to $\alpha \in [0,3]$, which is depicted as the violet region at the left-hand side of Fig.~\ref{fig:lyapunovDiag}. The bifurcation occurs near the midpoint of this region, where the Lyapunov exponent approaches~$0$. However, the isolating neighborhoods that separate the trajectories appear only when entering region (b), which corresponds to $\alpha = 3$.

\subsubsection{A period-$2$ limit cycle, two saddle fixed points and an unstable fixed point}
\label{sec:regionC}

\begin{figure}[htbp]
    \centering
    \begin{subfigure}[c]{0.48\linewidth}
        \centering
        \includegraphics[scale=1]{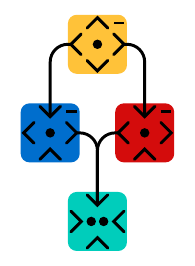}
    \end{subfigure}
    \begin{subfigure}[c]{0.48\linewidth}
        \centering
        \includegraphics[scale=0.5]{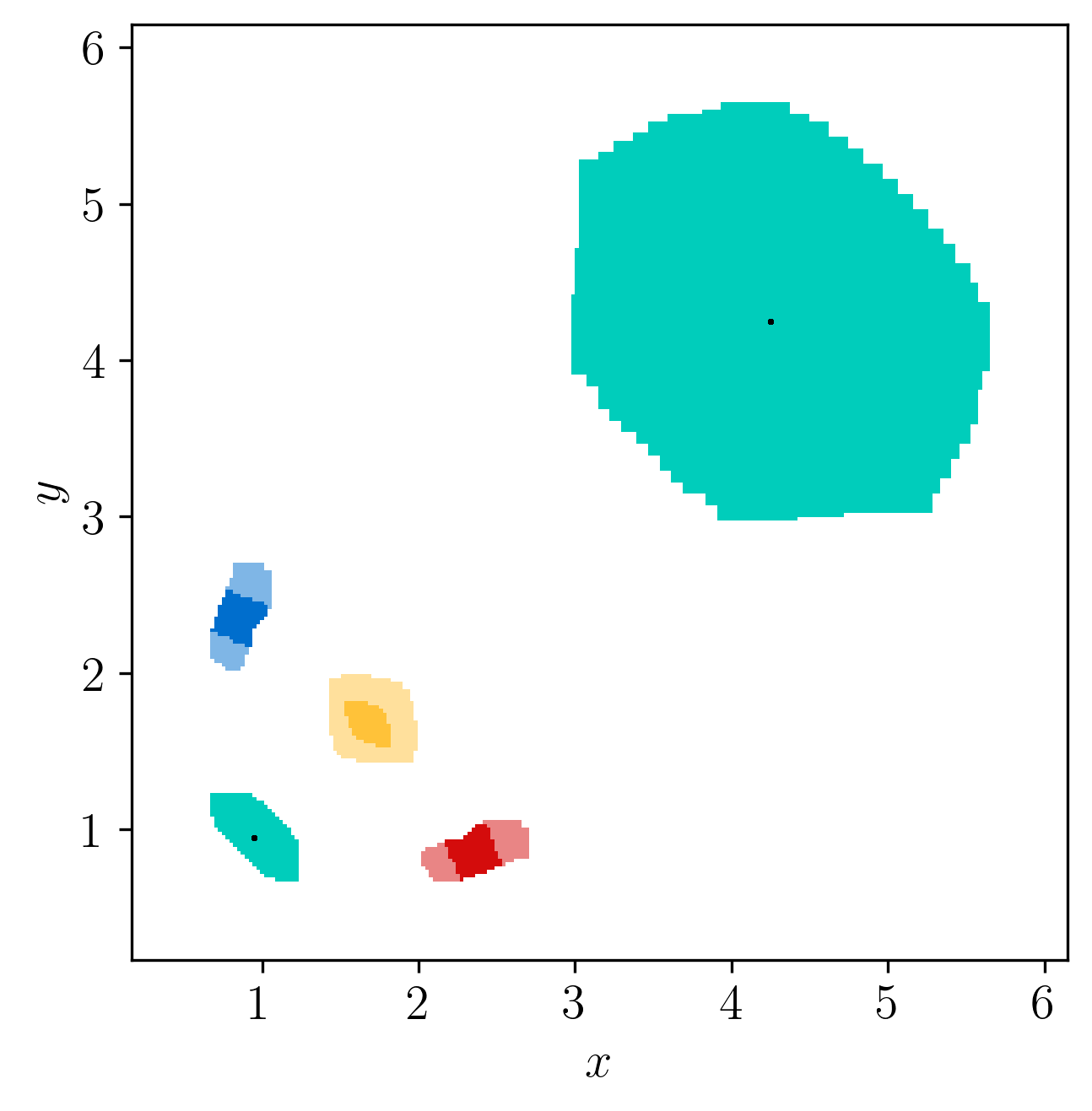}
    \end{subfigure}
    \caption{Conley--Morse graph and the isolating neighborhoods along with their exit sets found in the phase space for $\alpha_1, \alpha_2 \in [7, 7.5]$ in region (c).}
    \label{fig:15_15}
\end{figure}

A further increase in parameter values leads to region (c) containing the segment $\alpha_1 = \alpha_2 \in [5,16.5]$.
The large unstable numerical Morse set shown in Fig.~\ref{fig:8_8} splits into three sets illustrated in Fig.~\ref{fig:15_15}. This change corresponds to a pitchfork bifurcation.

The first set that appears at the top of the CM graph has two unstable directions and corresponds to an unstable fixed point (source). This information is provided by the Conley index and can also be deduced from the shape of the exit set (shown in pale yellow) that surrounds the isolating neighborhood entirely.
This set is accompanied by two numerical Morse sets (shown in blue and red), each with one unstable direction, corresponding to two saddle fixed points (the exit sets are shown in pale shades of blue and red, respectively).
Finally, a numerical Morse set consisting of two connected components is the largest set in this decomposition. It is an attractor in the sense that $f(M) \subset \interior M$ and thus the exit set is empty. The two connected components are mapped to each other. This set looks like an isolating neighborhood of a stable periodic orbit with period~$2$.

\llabel{pitchforkBegin}Looking at the shape and size of the unstable isolating neighborhood shown in yellow in Fig.~\ref{fig:8_8}, we might suspect that the pitchfork bifurcation took place at an earlier time, maybe even before the period doubling bifurcation. However, it has only been detected now with the CMAD method applied at fixed resolution.\llabel{pitchforkEnd}

\llabel{flipBegin}The flip in phase space caused by the periodic point of period $2$ is reflected in the three unstable fixed points (a source and two saddles).\llabel{flipEnd}


\subsubsection{Bistability in region (d)}
\label{sec:regionD}

\begin{figure}[htbp]
    \centering
    \begin{subfigure}[c]{0.48\linewidth}
        \centering
        \includegraphics[scale=1]{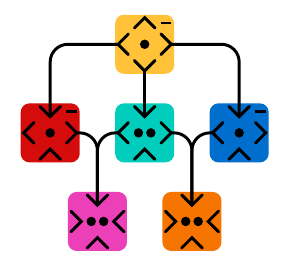}
    \end{subfigure}
    \begin{subfigure}[c]{0.48\linewidth}
        \centering
        \includegraphics[scale=0.5]{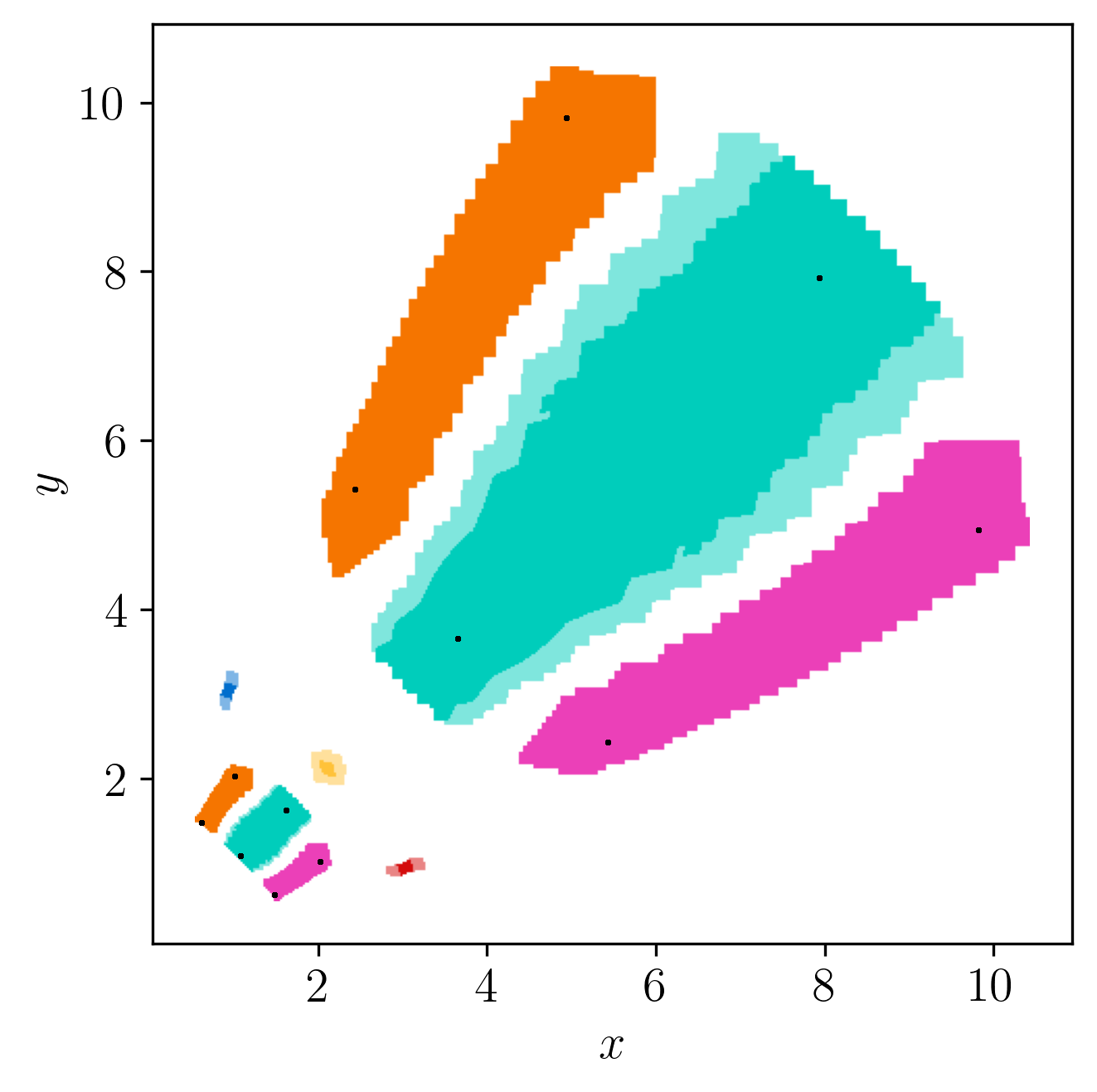}
    \end{subfigure}
    \caption{Conley--Morse graph and the isolating neighborhoods along with their exit sets found in the phase space for $\alpha_1,\alpha_2\in[17.5,18]$ in region (d). Attractors (also those restricted to the diagonal) found with numerical simulations for $\alpha_1=\alpha_2=17.5$ are indicated in black.}
    \label{fig:35_35}
\end{figure}

As $\alpha = \alpha_1 = \alpha_2$ increases within region (c), the attracting $2$-periodic isolating neighborhood (shown in cyan in Fig.~\ref{fig:15_15}) grows in size until it eventually splits into three numerical Morse sets in region (d).
The CM graph for $\alpha \in [16.5, 19]$ in region (d) reflects relatively complex dynamics of the system, as shown in Fig.~\ref{fig:35_35}.

\llabel{pitchfork2begin}The three sets (shown in cyan, orange, and pink) into which the previously stable $2$-periodic set was split look like the result of a pitchfork bifurcation of the stable periodic orbit of period $2$. Numerical simulations (described below) reveal an additional period doubling bifurcation, as the three products of the pitchfork bifurcation contain periodic orbits of period $4$, with two points in each connected component of the corresponding isolating neighborhoods. However, our method does not provide means to determine the order in which these two bifurcations occur.\llabel{pitchfork2end}

The isolating neighborhood of an unstable fixed point that was at the top of the CM graph in region (c) is still present in region (d), shown in yellow, accompanied by two small numerical Morse sets (isolating neighborhoods of saddle fixed points, shown in red and blue).\llabel{saddle} The map on these three sets reverses the orientation, which is related\llabel{revOrient} to the periodic points of period $2$ described above.



The numerical simulations for region (d) were performed for $\alpha_1=\alpha_2=17.5$ and the stable periodic orbits found are indicated with black dots inside the sets shown in orange and in pink. The attractors are very small in comparison to the isolating neighborhoods and have a different nature (periodic orbits of period $4$, not~$2$). At this point, the difference between our topological approach and the approach based on numerical simulations becomes apparent. Although simulations are easier to perform and display both attractors more accurately, isolating neighborhoods reveal features that would be difficult to determine otherwise:
The shape and size of the isolating neighborhoods reflect the strength of hyperbolicity of the isolated invariant sets.

As explained in the next paragraph, we used the symmetry in model \eqref{eq:model} and the fact that we are now considering $\alpha_1=\alpha_2$ and $\beta_1=\beta_2$ to conduct numerical simulations restricted to the diagonal that revealed an unstable period-$4$ orbit inside the large unstable numerical Morse set drawn in cyan in the center of Fig.~\ref{fig:35_35}. \llabel{orbit4begin}Note that the Morse set consists of $2$ connected components and the orbit consists of $4$ points, with $2$ points in each of the two components. One could say that this shows that our coarse computations are one period doubling bifurcation behind the results of numerical simulations in this case.\llabel{orbit4end}

Let us explain how we obtained the periodic orbit inside the unstable isolating neighborhood shown in cyan in Fig.~\ref{fig:35_35}. Recall that numerical simulations usually cannot portray unstable orbits. \llabel{numOrbitBegin}Even if we knew one point on such an orbit, rounding errors encountered when computing its image would drift the numerically computed orbit away from it.\llabel{numOrbitEnd} However, in our setting, model \eqref{eq:model} restricted to the diagonal $D=\{(x,y)\colon x=y\}$ becomes a one-dimensional system with dynamics contained fully within $D$, so we were computing numerical orbits for this constrained system to find the periodic orbit\llabel{numOrbitOK} that was stable there. We remark that such conditions allow certain unstable invariant sets to be visualized also in some of the following figures.

\subsubsection{Two stable orbits of period $4$ in region (e)}
\label{sec:regionE}

\begin{figure}[htbp]
    \centering
    \begin{subfigure}[b]{0.48\linewidth}
        \centering
        \includegraphics[scale=1]{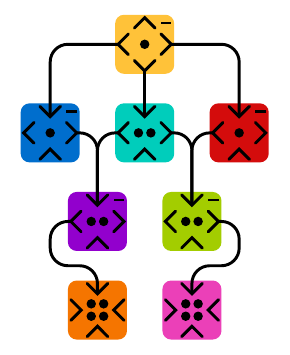}
    \end{subfigure}
    \begin{subfigure}[b]{0.48\linewidth}
        \centering
        \includegraphics[scale=0.5]{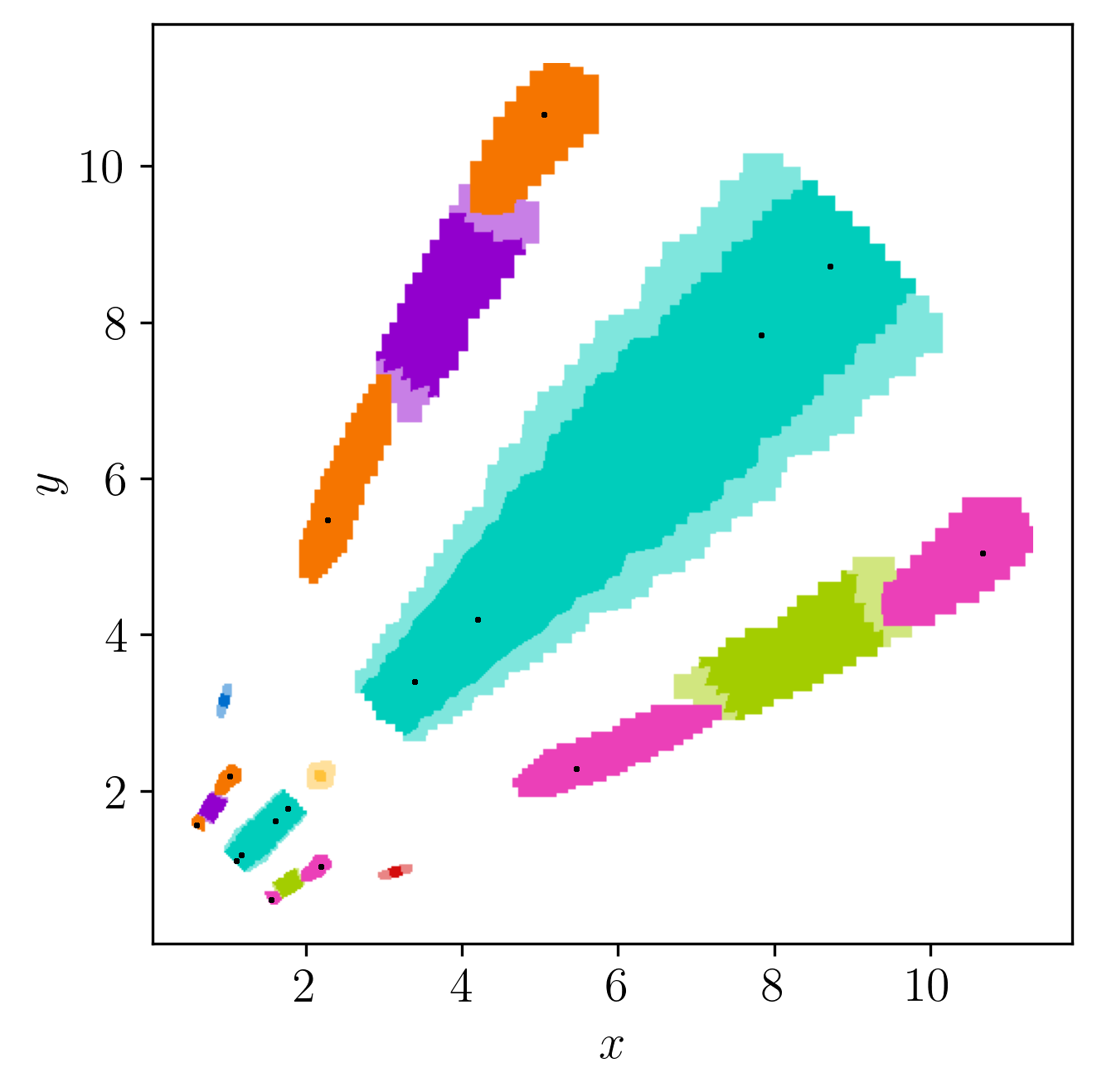}
    \end{subfigure}
    \caption{Conley--Morse graph and the isolating neighborhoods along with their exit sets found in the phase space for $\alpha_1,\alpha_2\in[20,20.5]$ in region (e).}
    \label{fig:40_40}
\end{figure}

Increasing the parameter $\alpha = \alpha_1 = \alpha_2$ further leads to region (e) containing the segment $\alpha_1 = \alpha_2 \in [19, 22.5]$. The new CM graph and a corresponding phase space portrait are shown in Fig.~\ref{fig:40_40}.

The main difference compared to region (d) is a stronger separation of the two attractors (one shown in orange, the other one in pink) from the large unstable numerical Morse set (shown in cyan) between them. \llabel{period4begin}Each of these two attractors now has $4$ connected components as a consequence of a period-doubling bifurcation of each of the two period-$2$ orbits (shown in purple and in green) that are now unstable.\llabel{period4end}

\subsubsection{Bistability and chaos in region (f)}
\label{sec:regionF}

\begin{figure}[htbp]
    \centering
    \begin{subfigure}[c]{0.48\linewidth}
        \centering
        \includegraphics[scale=1]{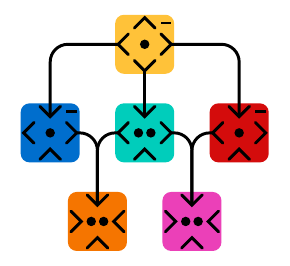}
    \end{subfigure}
    \begin{subfigure}[c]{0.48\linewidth}
        \centering
        \includegraphics[scale=0.5]{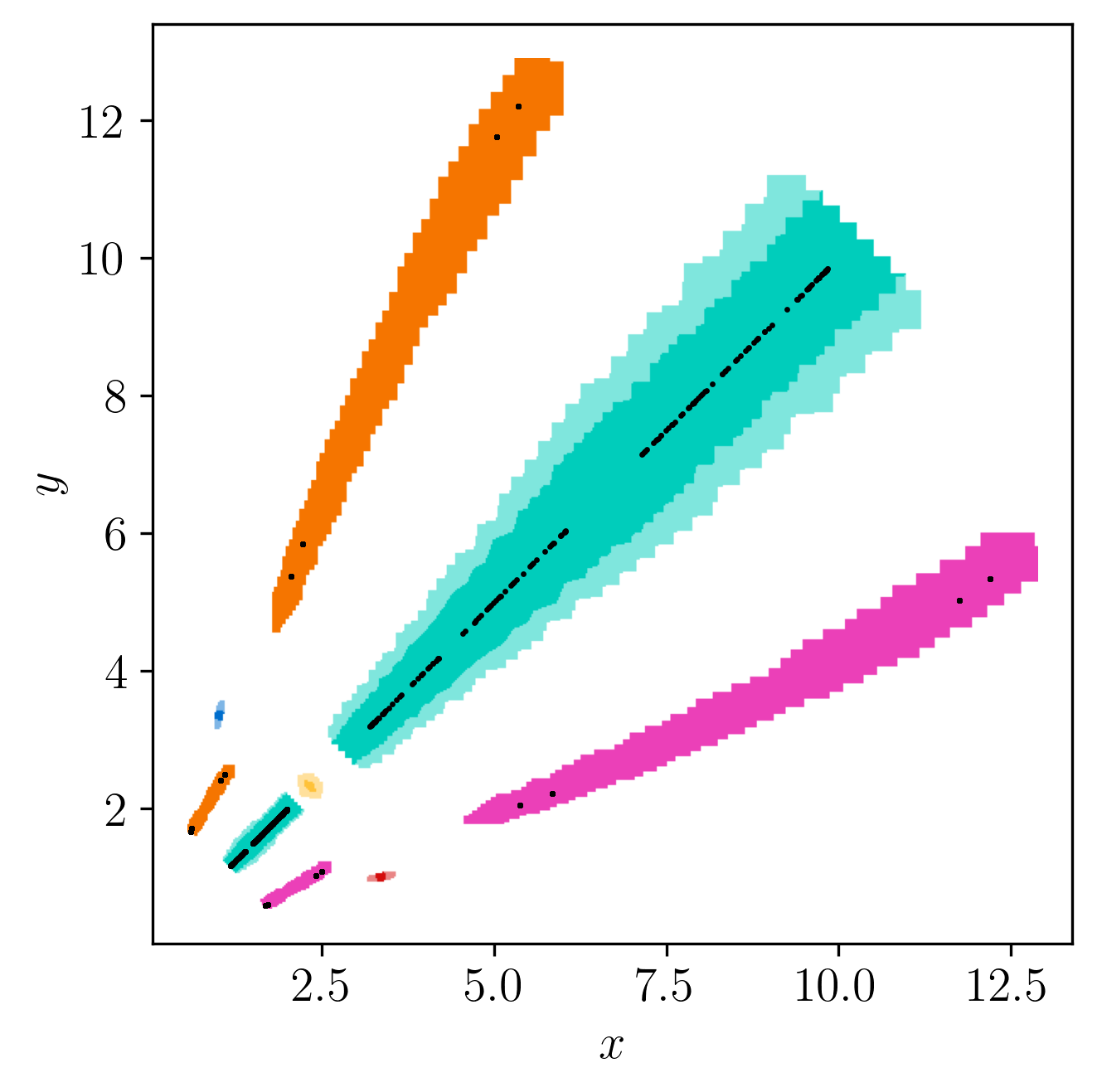}
    \end{subfigure}
    \caption{Conley--Morse graph and the isolating neighborhoods along with their exit sets found in the phase space for $\alpha_1,\alpha_2\in[25,25.5]$ in region (f).}
    \label{fig:50_50}
\end{figure}

As we follow the parameters $\alpha_1 = \alpha_2$ further and enter region (f) that includes the segment $\alpha_1 = \alpha_2 \in [22.5, 36]$, the dynamics observed at the chosen resolution reverts to that seen earlier in region (d), as illustrated in Fig.~\ref{fig:50_50}, except that the two attractors are thinner and farther apart. The isolating neighborhoods containing the period-$2$ saddle points have merged back into the neighborhoods of the attractors.

\llabel{regionFbifBegin}This is what we see at the level of isolating neighborhoods. The actual dynamics becomes more complex. Indeed, the bifurcation diagram and the corresponding changes in the largest Lyapunov exponent shown in Fig.~\ref{fig:lyapunovDiag}\llabel{refToBifDiag} clearly indicate that a period-doubling cascade is taking place when we traverse with $\alpha_1 = \alpha_2$ through region (f) and chaotic dynamics emerges, interrupted by some periodic windows. However, the CMAD method cannot see this due to the limited resolution of the phase space.\llabel{regionFbifEnd}

The impression of a simplified dynamical structure, due to a smaller CM graph, is caused by the fact that the actual dynamics becomes more and more complex within a limited area of the phase space. Further period-doubling bifurcations, weaker hyperbolicity of stable periodic orbits, and larger images of grid elements in the phase space due to larger values of parameters of the map and larger expansion have a detrimental effect on the ability to construct disjoint isolating neighborhoods at fixed finite resolution. Therefore, large isolating neighborhoods are created that contain several smaller invariant sets and finer details become hidden.

Increasing the resolution of computations may reveal more invariant sets. In our specific case, it turns out that conducting the computations whose results are shown in Fig.~\ref{fig:50_50} for region (f) at twice higher resolution in the phase space (with the subdivision depth, as defined in Sec.~\ref{sec:resol}, increased from $12$ to $13$) yields a very similar picture of the dynamics to what is shown in Fig.~\ref{fig:40_40} for region~(e). Inspection of the bifurcation diagram shown in Fig.~\ref{fig:lyapunovDiag} suggests that, if the accuracy of the computations is high enough, we should be able to see an isolating neighborhood of the two attractors consisting of four connected components for $\alpha$ in most of the interval $[22.5,36]$.

\llabel{denial}

\subsubsection{Chaotic dynamics in region (g)}
\label{sec:regionG}

\begin{figure}[htbp]
    \centering
    \begin{subfigure}[c]{0.4\linewidth}
        \centering
        \includegraphics[scale=1]{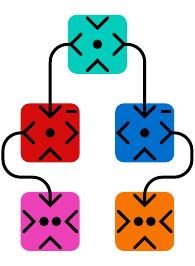}
    \end{subfigure}
    \hspace{0.05\linewidth}
    \begin{subfigure}[c]{0.5\linewidth}
        \centering
        \includegraphics[scale=0.5]{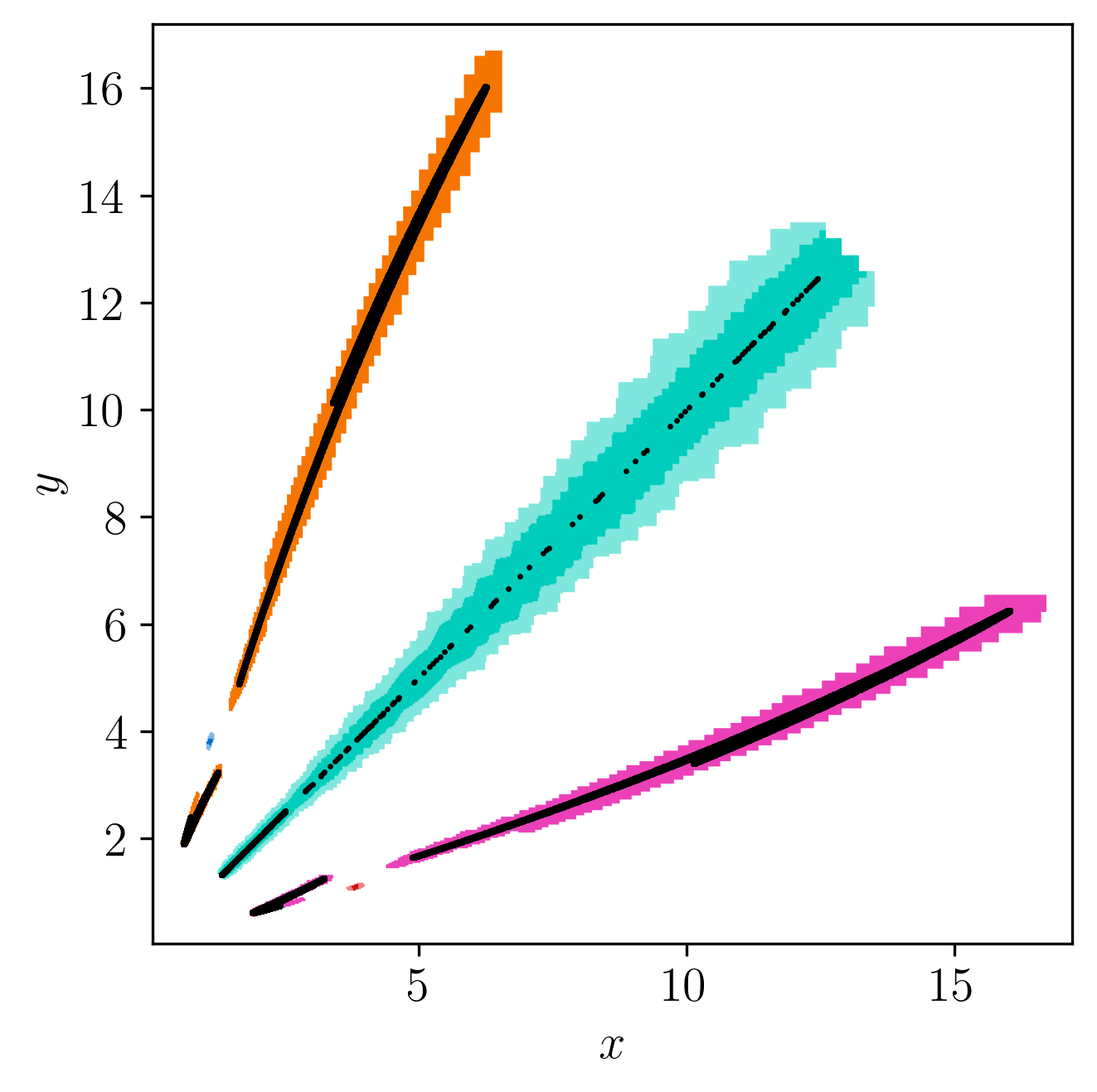}
    \end{subfigure}
    \caption{Conley--Morse graph and the isolating neighborhoods along with their exit sets found in the phase space for $\alpha_1, \alpha_2 \in [40, 40.5]$ in region (g).}
    \label{fig:80_80}
\end{figure}

As we continue increasing the parameters $\alpha_1 = \alpha_2$ and enter region (g) that contains the segment $\alpha_1 = \alpha_2 \in [37.5, 57.5]$, the CM graph becomes reduced further, as shown in Fig.~\ref{fig:80_80}.

The small neighborhood of the unstable fixed point (source) in region (f) shown in yellow in Fig.~\ref{fig:50_50} got merged with the large numerical Morse set shown in cyan into a large set with one stable and one unstable direction. This unstable fixed point may still be present inside the set, but at the chosen resolution it is not isolated anymore.

Each of the two $2$-periodic attractors is accompanied by a small numerical Morse set (a saddle fixed point; one shown in red, the other one in blue) between its two connected components.

\subsubsection{Further simplification of the CM graph in region (h)}
\label{sec:regionH}

\begin{figure}[htbp]
    \centering
    \begin{subfigure}[c]{0.4\linewidth}
        \centering
        \includegraphics[scale=1]{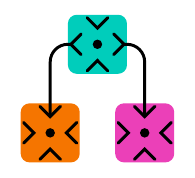}
    \end{subfigure}
    \hspace{0.05\linewidth}
    \begin{subfigure}[c]{0.5\linewidth}
        \centering
        \includegraphics[scale=0.5]{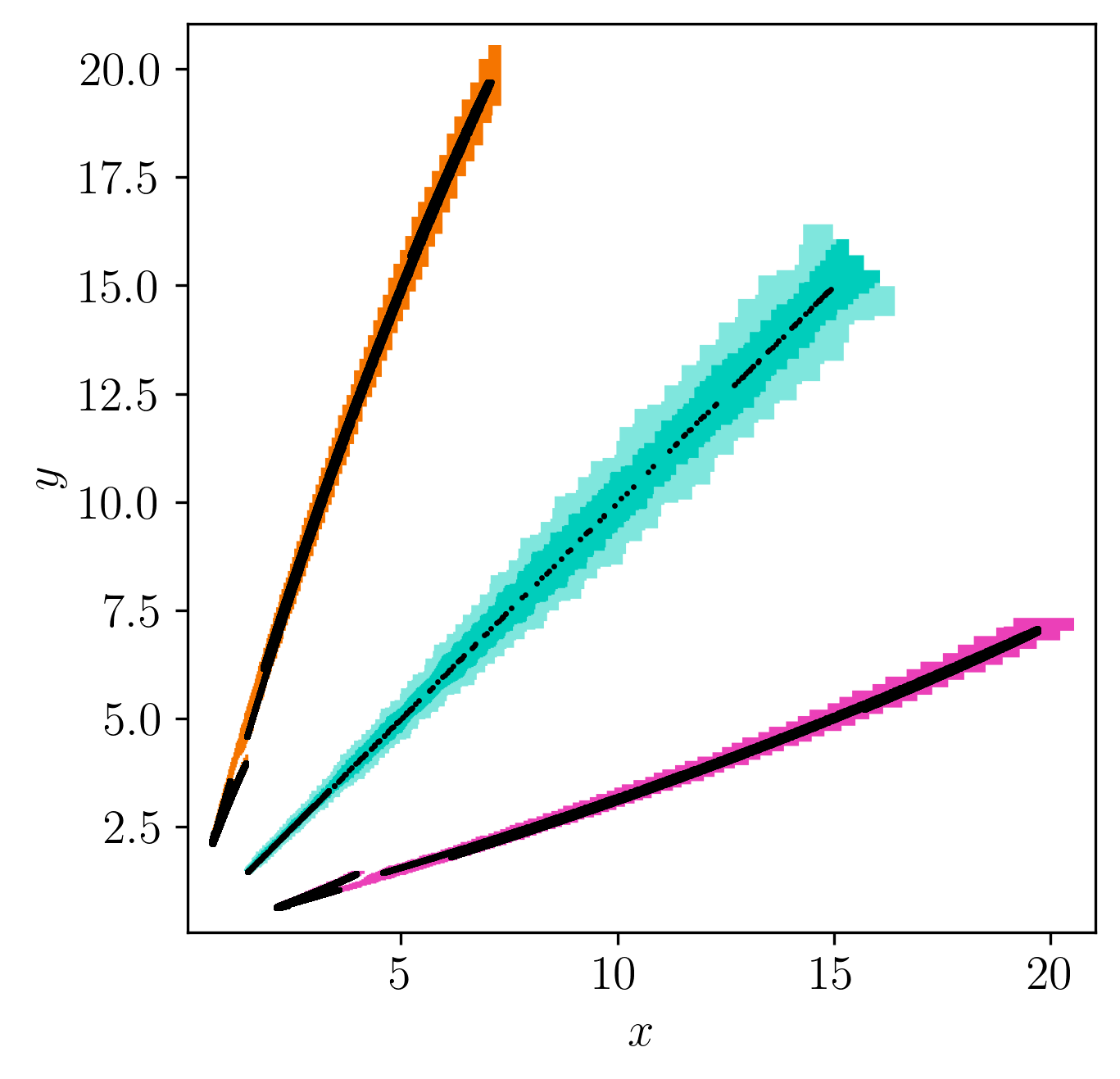}
    \end{subfigure}
    \caption{Conley--Morse graph and the isolating neighborhoods along with their exit sets found in the phase space for $\alpha_1, \alpha_2 \in [60, 60.5]$ in region (h).}
    \label{fig:120_120}
\end{figure}

Finally, in region (h) that contains the segment $\alpha_1 = \alpha_2 \in [57.5, 80]$, only three large numerical Morse sets were found.
The middle one with one unstable direction (shown in cyan) looks the same as in region (g).
The isolating neighborhoods of each of the two attractors\llabel{twoAttractors} (one shown in orange, the other one in pink) consist of one connected component each, and the small unstable numerical Morse sets (saddle fixed points) visible in region (g) have been swallowed by them.
The bifurcation diagram and the largest Lyapunov exponent shown in Fig.~\ref{fig:lyapunovDiag} indicate the presence of chaotic attractors, and the attractors found in numerical simulations show that the isolating neighborhoods encompass them relatively tight.

\subsection{Behavior of the system for different values of \texorpdfstring{$\alpha_1$}{alpha1} and \texorpdfstring{$\alpha_2$}{alpha2}}
\label{sec:different}

Following the variety of qualitatively different types of dynamics found in continuation classes along the diagonal, discussed in Sec.~\ref{sec:diag}, let us now examine the dynamics and bifurcations that can be observed for parameters apart from the diagonal.
We will check how the system's dynamics changes as we approach the diagonal.

\subsubsection{Two period-doubling bifurcations}
\label{sec:region1}

To begin with, let us assume that $\alpha_1$ is approximately twice as large as $\alpha_2$.  
In this case, we pass through regions (a), (b), (1), and (2) in the continuation diagram shown in Fig.~\ref{fig:contDiag}.

\begin{figure}[htbp]
    \centering
    \begin{subfigure}[b]{0.22\linewidth}
        \centering
        \includegraphics[scale=1]{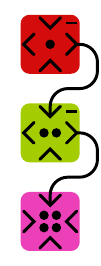}
    \end{subfigure}
    \hspace{0.05\linewidth}
    \begin{subfigure}[b]{0.5\linewidth}
        \centering
        \includegraphics[scale=0.5]{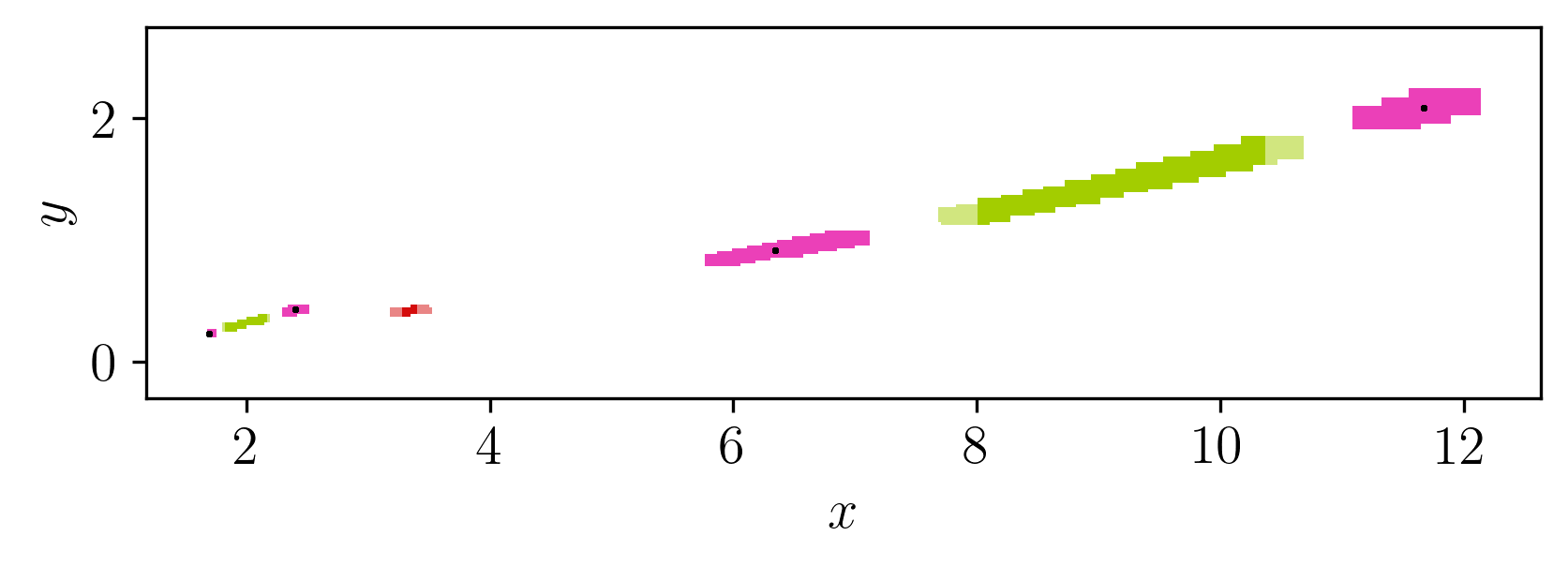}
    \end{subfigure}
    \caption{Conley--Morse graph and the isolating neighborhoods along with their exit sets found in the phase space for $\alpha_1 \in [22.5, 23]$, $\alpha_2 \in [10, 10.5]$ in region (1).}
    \label{fig:45_20}
\end{figure}

Recall that for $(\alpha_1,\alpha_2)$ moving from region (a) to region (b), we observed a period-doubling bifurcation (see Sec.~\ref{sec:regionB}) and the numerical Morse decomposition in region (b) suggests a stable periodic orbit of period $2$ and an unstable saddle fixed point, shown in Fig.~\ref{fig:8_8}. As the parameters $\alpha_1$ and $\alpha_2$ increase, these components become increasingly elongated, until another period-doubling bifurcation occurs and we see isolating neighborhoods of a saddle fixed point, a saddle periodic orbit of period $2$, and a stable periodic orbit of period $4$. The dynamics in this continuation class, region (1), is shown in Fig.~\ref{fig:45_20} for sample parameters chosen from this region.

\subsubsection{A coarse Morse decomposition}
\label{region2}

\begin{figure}[htbp]
    \centering
    \begin{subfigure}[b]{0.15\linewidth}
        \centering
        \includegraphics[scale=1]{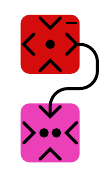}
    \end{subfigure}
    \hspace{0.05\linewidth}
    \begin{subfigure}[b]{0.6\linewidth}
        \centering
        \includegraphics[scale=0.5]{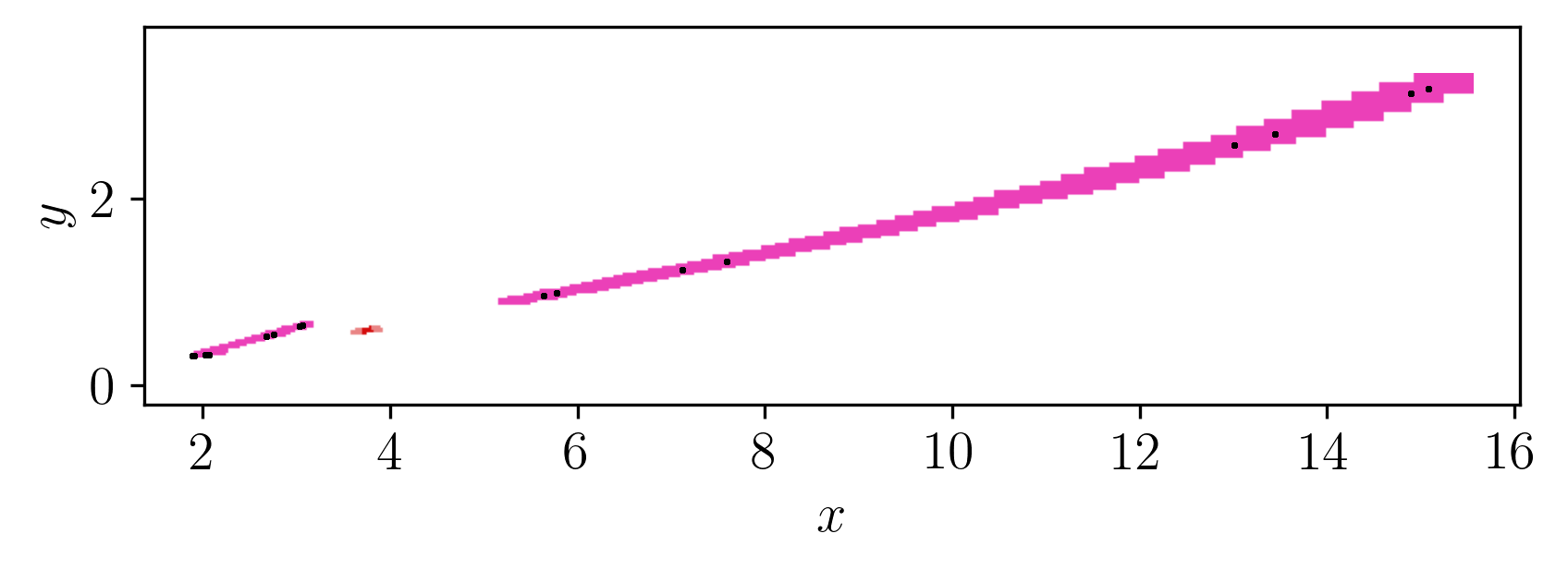}
    \end{subfigure}
    \caption{Conley--Morse graph and the isolating neighborhoods along with their exit sets found in the phase space for $\alpha_1\in[35,35.5]$, $\alpha_2\in[20,20.5]$ in region (2).}
    \label{fig:70_40}
\end{figure}

Increasing the parameters $\alpha_1$ and $\alpha_2$ further leads to region (2), where the CM graph becomes reduced, as shown in Fig.~\ref{fig:70_40}, but numerical simulations show evidence of more period-doubling bifurcations; in the figure, one can see a stable periodic orbit of period $16$, as there are $8$ black points in each of the two components of the isolating neighborhood shown in pink.

Note that bistability does not occur in regions (1) and (2).  
A significant difference between the values of $\alpha_1$ and~$\alpha_2$ causes only one attractor to appear in the system at a time.
The second attractor will begin to form in a series of bifurcations as we approach the diagonal, which will be illustrated in the following analysis, passing through regions (3), (4), and (5).

Since $\alpha_1$ and $\alpha_2$ represent the expression speeds of the two genes, the lack of bistability can be interpreted from a biological point of view as a permanent domination of one gene over the other due to the higher production of its corresponding protein.
Reducing the difference between the parameters (and thus approaching the diagonal in the continuation diagram) results in the gradual appearance of more complex structures in regions in which both genes are able to be expressed at comparable levels.

\subsubsection{A numerical Morse set with the trivial Conley index}
\label{sec:region3}

\begin{figure}[htbp]
    \centering
    \begin{subfigure}[b]{0.20\linewidth}
        \centering
        \includegraphics[scale=1]{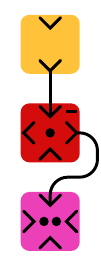}
    \end{subfigure}
    \hspace{0.05\linewidth}
    \begin{subfigure}[b]{0.6\linewidth}
        \centering
        \includegraphics[scale=0.5]{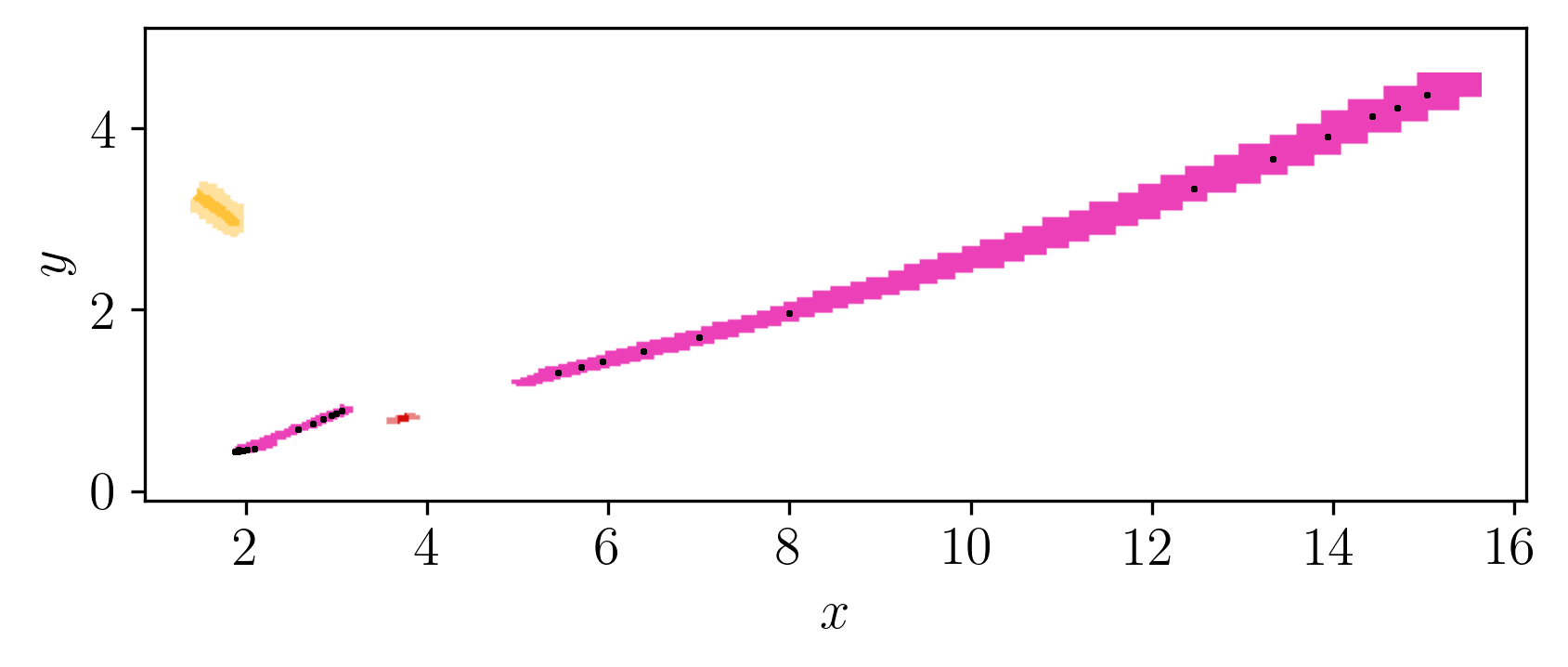}
    \end{subfigure}
    \caption{Conley--Morse graph and the isolating neighborhoods along with their exit sets found in the phase space for $\alpha_1 \in [35, 35.5]$, $\alpha_2 \in [27, 27.5]$ in region (3).}
    \label{fig:70_54}
\end{figure}

When moving from region (2) to region (3), a new numerical Morse set appears at some distance from the previous numerical Morse sets, as shown in Fig.~\ref{fig:70_54}.  
Its Conley index is trivial, so its invariant part may be empty.
This means that trajectories may simply pass through this set.
However, inclusion of this isolating neighborhood in the constructed CM graph indicates a local slow-down in the dynamics. This typically happens before some types of bifurcation take place, as we show in Fig.~\ref{fig:sn_bifur}, discussed in Sec.~\ref{sec:bifurcations}. Indeed, after a slight increase of $\alpha_2$ this neighborhood splits into two sets, as described in Sec.~\ref{sec:region4}.

\subsubsection{A saddle--node bifurcation}
\label{sec:region4}

\begin{figure}[htbp]
    \centering
    \begin{subfigure}[b]{0.4\linewidth}
        \centering
        \includegraphics[scale=1]{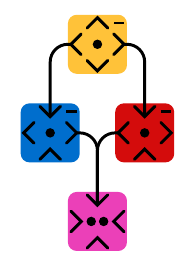}
    \end{subfigure}
    \hspace{0.05\linewidth}
    \begin{subfigure}[b]{0.5\linewidth}
        \centering
        \includegraphics[scale=0.5]{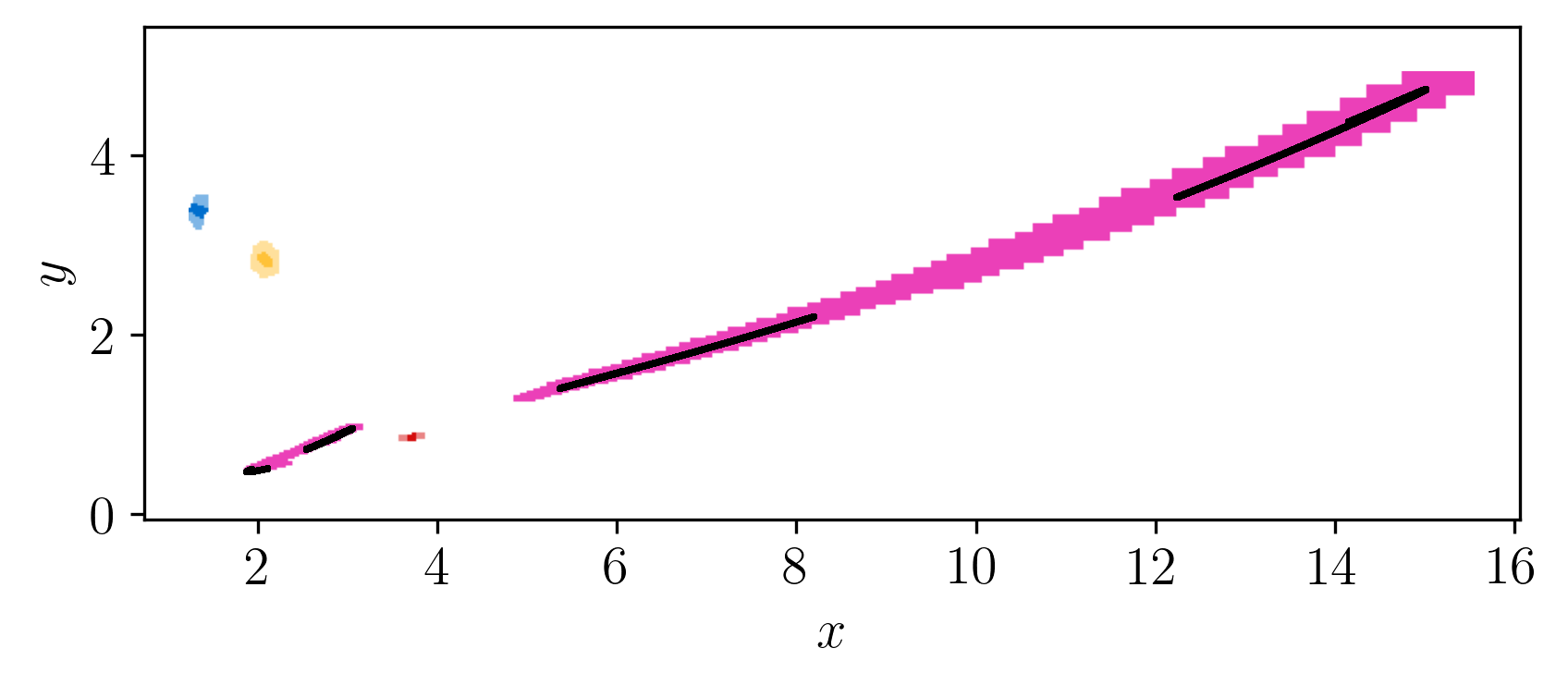}
    \end{subfigure}
    \caption{Conley--Morse graph and the isolating neighborhoods along with their exit sets found in the phase space for $\alpha_1 \in [35, 35.5]$, $\alpha_2 \in [29, 29.5]$ in region (4).}
    \label{fig:70_58}
\end{figure}

Figure~\ref{fig:70_58} shows the CM graph and the phase space for a sample square of parameters $\alpha_1$ and $\alpha_2$ in region (4).
The previously observed numerical Morse set with the trivial Conley index gave rise to two numerical Morse sets with nontrivial indices: a saddle fixed point (shown in blue) and an unstable fixed point (a source; shown in yellow).
The dynamics resembles the case already observed in region (c), shown in Fig.~\ref{fig:15_15}, except that the currently observed stable isolating neighborhood (shown in pink) is thin and located towards bottom and right, as if it was only a subset of the thick period-$2$ attracting isolaging neighborhood shown in cyan in Fig.~\ref{fig:15_15}.

\subsubsection{Another saddle--node bifurcation}
\label{sec:region5}

\begin{figure}[htbp]
    \centering
    \begin{subfigure}[b]{0.35\linewidth}
        \centering
        \includegraphics[scale=1]{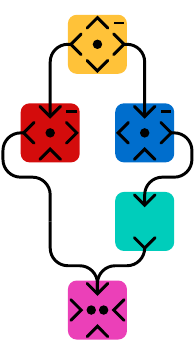}
    \end{subfigure}
    \hspace{0.05\linewidth}
    \begin{subfigure}[b]{0.55\linewidth}
        \centering
        \includegraphics[scale=0.5]{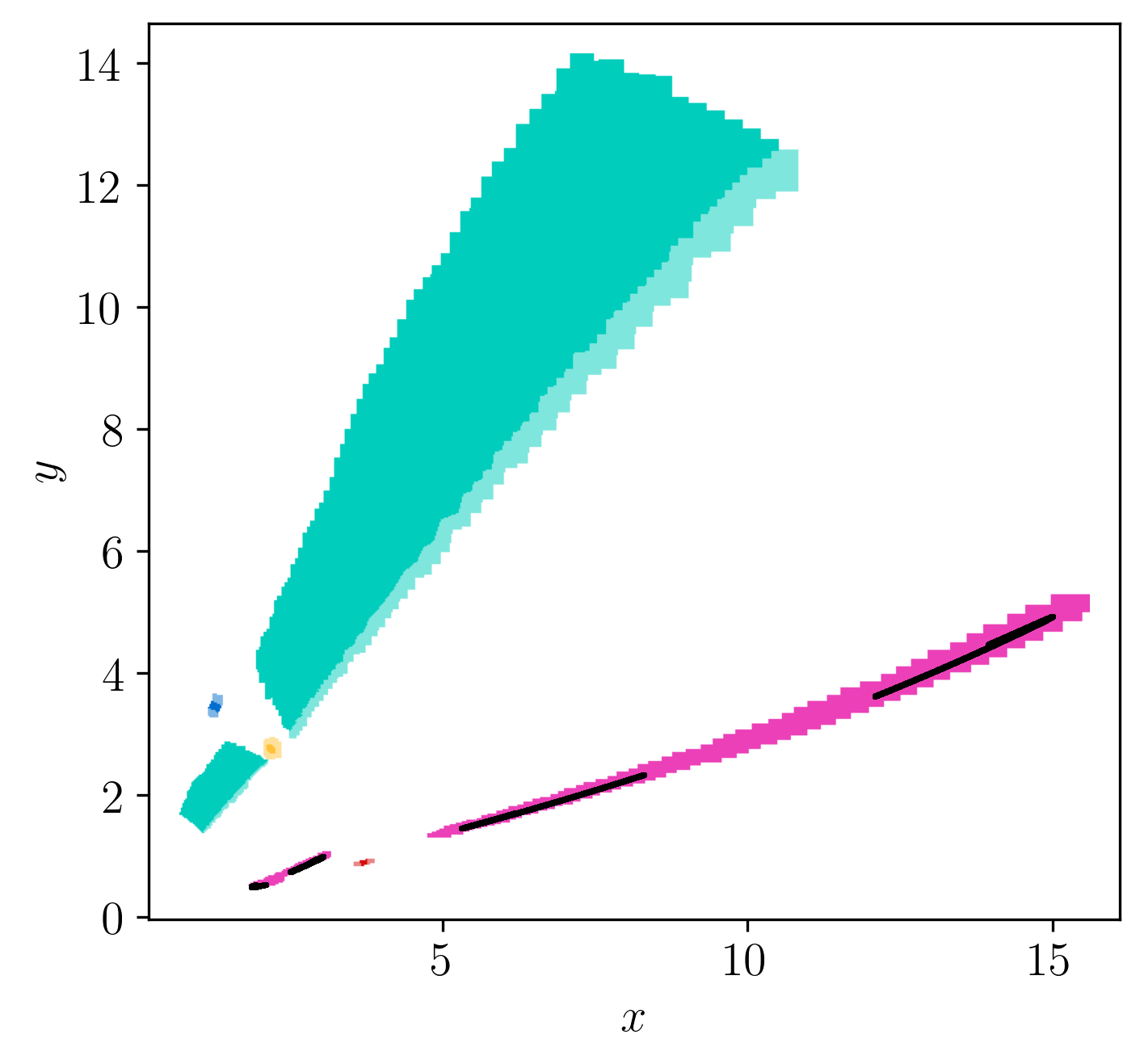}
    \end{subfigure}
    \caption{Conley--Morse graph and the isolating neighborhoods along with their exit sets found in the phase space for $\alpha_1 \in [35, 35.5]$, $\alpha_2 \in [30, 30.5]$ in region (5).}
    \label{fig:70_60}
\end{figure}

Further approaching the diagonal of the continuation diagram with the parameters $\alpha_1$ and $\alpha_2$ leads to region (5).
In this region, a new large numerical Morse set appears with the trivial Conley index, as shown in Fig.~\ref{fig:70_60}.
Judging by the exit set located entirely at the right--bottom edge of the set, the trajectories seem to merely pass through this isolating neighborhood, but the slowdown is significant enough for the CMAD method to treat it as a separate numerical Morse set.

Further approaching the diagonal leads to region (f), in which this large neighborhood is split into a second attractor (the node) accompanied by a numerical Morse set with one unstable direction (the saddle), both sets corresponding to periodic orbits of period $2$ in our case.

It is interesting to see that we have observed two paths leading to bistability in the Andrecut--Kauffman system obtained by different sequences of bifurcations from the most simple dynamics observed in region (a), corresponding to a stable fixed point. The first path, following the diagonal and described in Sec.~\ref{sec:diag}, consists mainly of splitting the existing numerical Morse sets into smaller ones, while the second path, leading at some distance from the diagonal and approaching it from a side, relies on the emergence of new pairs of numerical Morse sets through saddle--node bifurcations.

\subsection{Impact of the phase space resolution on the observed dynamics}
\label{sec:resol}

The resolution imposed by the grid taken in the phase space is a parameter that, combined with the size of the rectangular area $B$, has the most pronounced influence on the results provided by the CMAD method. In our case, the resolution is determined by the number $d \in \mathbb{N}$ of subdivisions of $B$ in each direction conducted during the gradual refinement approach. We shall call this number a \emph{subdivision depth}.

If we choose a number $d \in \mathbb{N}$ as the subdivision depth of $B$ then the set $B$ will be effectively subdivided into $2^d$ intervals in each direction. In our case of the two-dimensional model, the resulting grid will consist of $4^d$ elements.

The choice of a subdivision depth and thus the resolution in the phase space directly affects the running time of the computations, memory usage, and the obtained results. The higher the resolution, the more accurate the map representation is computed, and thus the richer the description of the dynamics (thanks to a more complex CM graph and tighter isolating neighborhoods). However, this also implies considerably higher cost of the computations, both in terms of time and memory, especially if large numerical Morse sets are constructed (that is, built of a large number of grid elements).

Let us now analyze an example that illustrates the effect of a chosen subdivision depth on the qualitative properties of the dynamical system that we observe.

Consider model \eqref{eq:model} with $\alpha_1,\alpha_2\in[6,6.5]$, $\varepsilon\in[0.94,0.95]$, $\beta_1=\beta_2=0.2$ and $n=3$.
Take $B = [0,101] \times [0,101]$.
In what follows, we show the resulting CM graphs and the corresponding numerical Morse sets constructed for three different values of the subdivision depth: $d \in \{8, 10, 11\}$.

\begin{figure}[htbp]
    \centering
    \begin{subfigure}[c]{0.2\linewidth}
        \centering
        \includegraphics[scale=1]{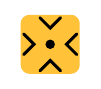}
    \end{subfigure}
    \hspace{0.05\linewidth}
    \begin{subfigure}[c]{0.5\linewidth}
        \centering
        \includegraphics[scale=0.5]{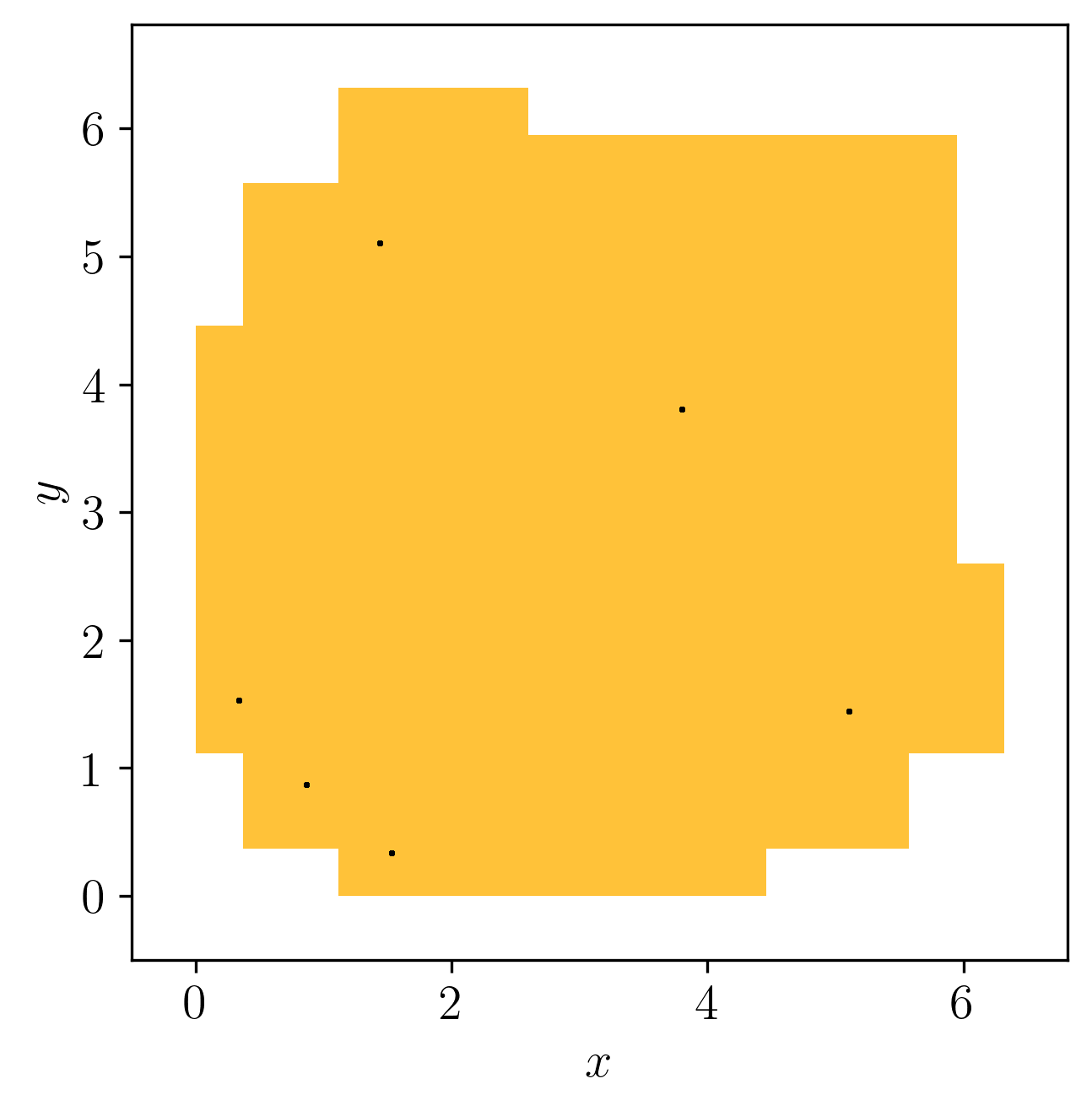}
    \end{subfigure}
    \caption{Conley--Morse graph and the isolating neighborhood found in the example discussed in Sec.~\ref{sec:resol} with subdivision depth $d=8$.}
    \label{fig:d8}
\end{figure}

Figure~\ref{fig:d8} shows a CM graph consisting of a single vertex and the corresponding numerical Morse set $M$. The set $B$ was subdivided into $2^8$ equal segments in both directions. It was proved that $f(M) \subset \interior M$. This set looks like an isolating neighborhood of a globally stable fixed point with relatively weak attraction.

\begin{figure}[htbp]
    \centering
    \begin{subfigure}[c]{0.18\linewidth}
        \centering
        \includegraphics[scale=1]{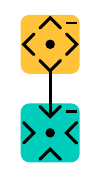}
    \end{subfigure}
    \hspace{0.05\linewidth}
    \begin{subfigure}[c]{0.5\linewidth}
        \centering
        \includegraphics[scale=0.5]{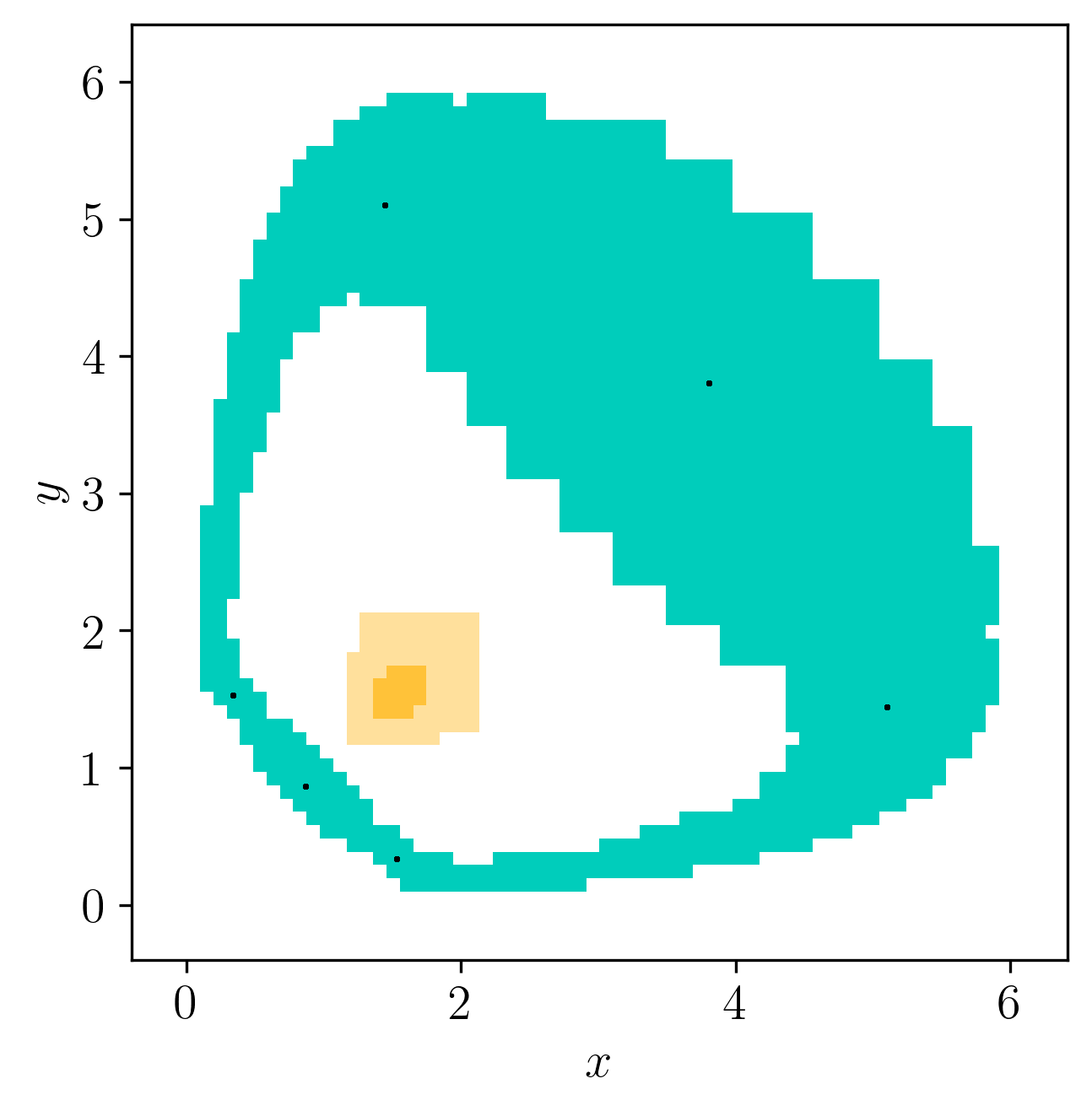}
    \end{subfigure}
    \caption{Conley--Morse graph and the isolating neighborhoods found in the example discussed in Sec.~\ref{sec:resol} with subdivision depth $d=10$.}
    \label{fig:d10}
\end{figure}

Figure~\ref{fig:d10} shows a CM graph consisting of two vertices and the two corresponding numerical Morse sets $M_1$ and $M_2$.
The set $B$ was subdivided into $2^{10}$ equal segments in both directions in this case, which implies a fourfold higher resolution in each direction or a 16-fold increase in the number of grid elements into which $B$ was effectively subdivided.

The set $M_1$ in the middle is a repeller with two unstable directions (an unstable fixed point), surrounded by an attracting ring (or circle) $M_2$. It was proved that $f(M_2) \subset \interior M_2$.

In this way, an increase in resolution allowed us to see a more detailed picture of the dynamics. In particular, we know that solutions tend to leave the middle part of the isolating neighborhood $M$ found previously and then stay in the ring. The Conley index provides somewhat puzzling information that the dynamics on the ring reverses the orientation of the space, so the map looks like a reflection of the ring about a line. This rules out the possibility of solutions running around the ring in circular motion, so it is not an invariant circle. What is the actual dynamics then?

\begin{figure}[htbp]
    \centering
    \begin{subfigure}[c]{0.35\linewidth}
        \centering
        \includegraphics[scale=1]{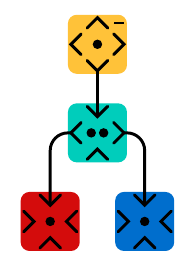}
    \end{subfigure}
    \hspace{0.05\linewidth}
    \begin{subfigure}[c]{0.5\linewidth}
        \centering
        \includegraphics[scale=0.5]{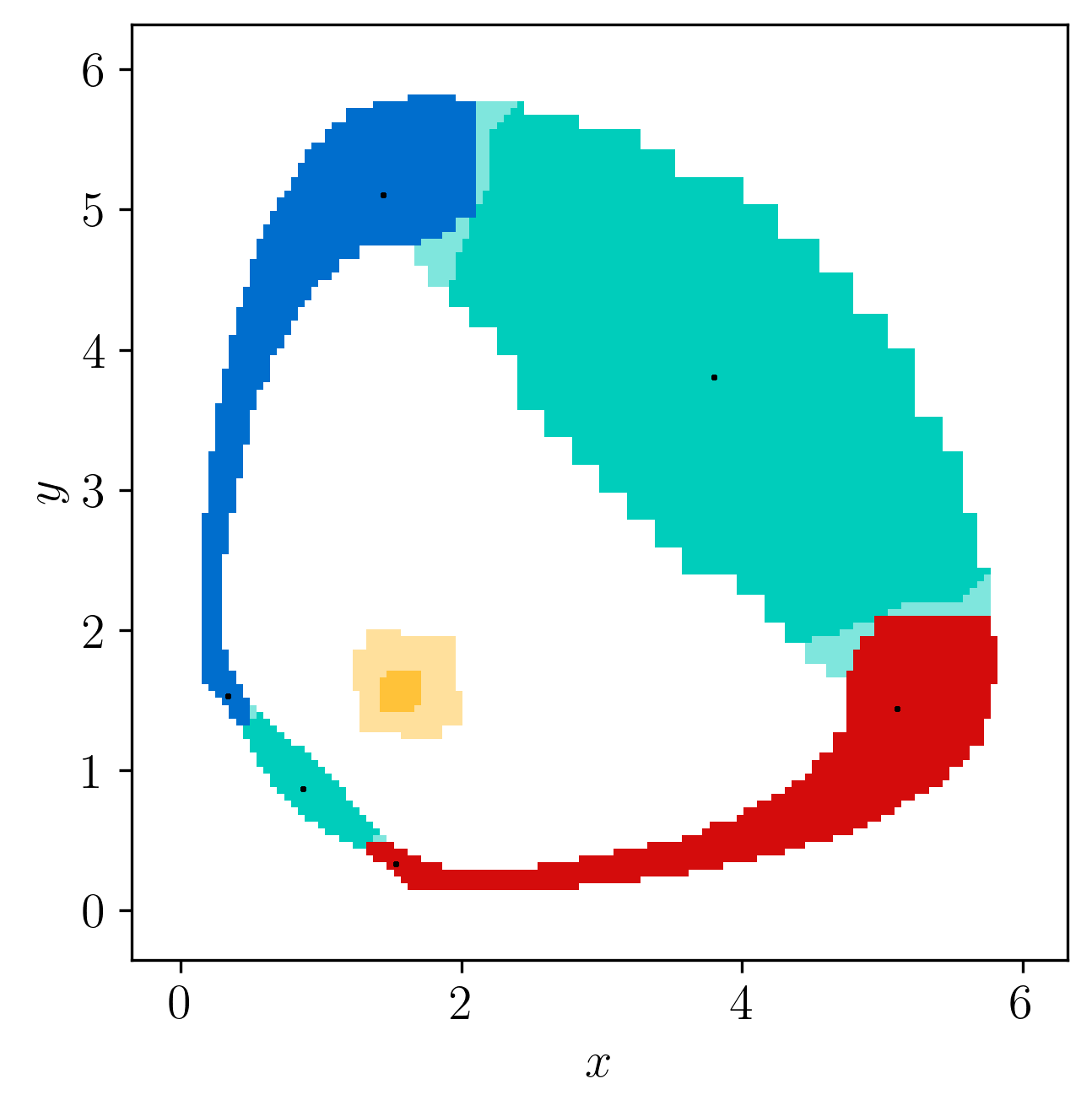}
    \end{subfigure}
    \caption{Conley--Morse graph and the isolating neighborhoods found in the example discussed in Sec.~\ref{sec:resol} with subdivision depth $d=11$.}
    \label{fig:d11}
\end{figure}

Figure~\ref{fig:d11} shows a CM graph consisting of four vertices and the corresponding four numerical Morse sets.  
The set $B$ was subdivided into $2^{11}$ equal segments in both directions in this case.

In addition to the previously seen isolating neighborhood $M_1$ (shown in yellow) of an unstable fixed point in the middle, we now see a clear picture of the dynamics that explains the reason for the orientation change. Namely, the numerical Morse set $M_2$ (shown in cyan) is an isolating neighborhood of a saddle periodic orbit of period $2$. It consists of two connected components, each mapped to the other. The remaining part of the ring is nearly completely filled with two attractors $M_3$ (shown in red) and $M_4$ (shown in blue). In particular, the system exhibits bistability.

This example shows that the subdivision depth of the phase space may have a significant impact on the results, and thus on our perception of the dynamics.
The choice of final resolution of the computations affects which structures can be isolated with the use of the graph representation of the map on a particular grid.
In our case, we began with the suspicion of the existence of an attracting fixed point.  
Ultimately, it turned out that the initial understanding of the dynamics was replaced with a deeper insight into the internal structure of the initially found isolating neighborhood.
In this way, we achieved a better understanding of the dynamics, although still at limited resolution.

This example also warns that the interpretation of the results should be carried out very cautiously to avoid premature conclusions.
One should keep in mind that the results are always based on computations conducted at a fixed finite resolution that include overestimates. In particular, we are not able to observe anything happening at a smaller scale than the assumed finite resolution.

Finally, we would like to emphasize the importance of the information on whether the map preserves orientation or not. Without such information, we could na\"ively suspect that the loop-shaped numerical Morse set observed in Fig.~\ref{fig:d10} could be a neighborhood of a stable circle with quasi-periodic motion on it. However, with a different choice of parameters, it is possible for the Andrecut--Kauffman model to exhibit behavior in which a ring-shaped attractor cannot be split like in our example, and the trajectories indeed run along the circle. Such a situation was originally found in \cite{subramani-2023}. An approximate attractor superimposed onto high resolution image of isolating neighborhoods for these parameter values is shown in Fig.~\ref{fig:subramani}.

\subsection{Importance of analyzing the Andrecut-Kauffman model}
\label{sec:importance}

Understanding the mechanisms and structural principles underlying gene expression regulation is fundamental to understanding the biological and chemical processes governed by genetic control. This knowledge is also essential for identifying the origins of various diseases and for developing effective therapeutic strategies. Disruptions in gene regulation can have profound pathological consequences.\llabel{test} For example, insufficient protein expression can cause premature cell differentiation and ultimately impair tissue development. On the other hand, protein overexpression is found in the pathogenesis of numerous diseases, including lung, ovarian, and colorectal cancers, as well as germ cell tumors~\cite{li-cao-li-jin-2018,liu-dai-du-2015}.

In particular, elevated expression levels can activate oncogenic signaling pathways, as observed in certain forms of leukemia, thus promoting tumor cell proliferation. Furthermore, dysregulated gene expression can contribute to metastasis by increasing the proliferative and migratory capacity of cancer cells. In breast cancer, for example, lower levels of expression of specific proteins have been associated with improved patient survival compared to cases exhibiting high expression~\cite{li-cao-li-jin-2018}. These findings emphasize the critical role of gene expression regulation not only in tumor initiation and progression but also in metastatic potential~\cite{liu-dai-du-2015}.

Given these observations, impaired regulation of gene expression emerges as a compelling target for therapeutic intervention. Overexpression and underexpression can both destabilize gene regulatory networks, potentially leading to severe pathological outcomes. Various mathematical models have been proposed to describe gene expression: the quasi-stationary approximation has been applied to self-regulating protein systems \cite{bartlomiejczyk-bodnar-2020}, while the effects of time delays on Hopf bifurcations and oscillations have been investigated in models with protein dimers and Hes1 \cite{bartlomiejczyk-bodnar-2023}. In addition, extended regulatory motifs such as the generalized p53–Mdm2 protein interaction network have been analyzed within a mathematical framework to investigate the stability and dynamic responses of gene expression systems \cite{piotorwska-bartlomiejczyk-bodnar-2018}.  These diverse approaches highlight that gene expression dynamics can be studied using both continuous models and models with time delays, providing a comprehensive understanding of regulatory system behavior.

Within this context, mathematical modeling of gene expression, such as the discrete-time two-gene Andrecut--Kauffman model, offers a valuable insight into the dynamic behavior of regulatory systems, facilitating a deeper understanding of their role in health and disease.

The discrete two-gene Andrecut--Kauffman model~\cite{andrecut_main}, based on the classical Boolean networks introduced by Kauffman in 1969~\cite{KAUFFMAN1969437}, provides a simplified yet powerful framework for studying gene regulation dynamics in systems with a limited number of components. It features threshold functions that describe gene interactions and is capable of modeling a wide range of dynamic behaviors, from stable fixed points through oscillations to deterministic chaos, \cite{Andrecut_2005, Andrecut_2006}.

Previous studies have mainly focused on local dynamical properties such as fixed points or limit cycles, yet modern analysis requires a global perspective on system dynamics in the phase space. The approach we propose allows us to capture complex relationships among trajectories and their topological organization, including unstable invariant sets, which is crucial for a comprehensive understanding of the behavior of a gene regulatory network.

Our results and the methods we propose may contribute in the future to a deeper understanding of the dynamics of gene regulatory networks and offer potential insights into biological processes such as cellular phenotype switching and tumor progression.

\section{Conclusion and final remarks}
\label{sec:final}

We have conducted a topological--numerical analysis of the dynamics in the Andrecut--Kauffman model \eqref{eq:model} at some fine but limited resolution fixed in the phase space. In this analysis, we shifted our focus away from periodic and chaotic dynamics in favor of a broader, more comprehensive view of the entire system. A description based on Conley--Morse graphs and continuation classes can be more intuitive and human--readable than raw data obtained from simulations.

\llabel{unstableBegin}It is also important to emphasize the fact that the approach we applied allows one to capture both the stable invariant sets that can be easily find in simulations, and also unstable invariant sets that are usually missed when using classical methods.
Although the unstable sets are not long-term states in which the system stabilizes, they are also of practical importance.
Namely, trajectories in close proximity of unstable sets, while short-lived on their own, can oftentimes be stabilized by gentle external influence to stay close to the unstable set for a prolonged time. Therefore, from the point of view of applications, they might be regarded almost as important as stable sets of the biological system, even though they require application of some control, for example, in the form of varying doses of a medication, to keep the system stable.\llabel{unstableEnd}

Theoretical rigor is one of the key features of the method that we applied. It provides a way of proving formal statements about dynamical systems that is especially useful if the complexity of the systems makes purely analytic methods infeasible. We feel that it might be worth applying this method to the analysis of other systems that have been analyzed only by means of numerical simulations so far. The finite--resolution results might confirm the credibility of the results of simulations and extend them to obtain a more comprehensive view of the dynamics.

Due to the limited resolution, we missed all the features that require a finer scale to be observed. However, our results are robust, and this makes them resistant to small perturbations. Indeed, every mathematical model describes a real biological system with limited precision, and there might also be some fluctuations or disturbances in the real system that are not included in the model, so any phenomena found in the model at a very small scale may not actually correspond to what is relevant to the biological system. Therefore, we feel that the finite-resolution approach to the analysis of a dynamical system that we applied might be a good choice to obtain meaningful results than using sophisticated analytic or numerical methods. See \cite{luzzatto-2011} for an in-depth discussion of this question and a possible way to formalize the notion of finite resolution dynamics.

\llabel{combining}We would like to note that the CMAD method does not \emph{replace} numerical simulations but provides a different class of results that can be complemented by numerical simulations. Such an integrated approach emphasizes the value of combining topological methods with numerical simulations in studying applied dynamical systems.

Despite extensive analysis in this and several previous works, the majority of the parameter space for the Andrecut--Kauffman model remains an unexplored territory. The case of $\beta_1\neq\beta_2$ appears exceptionally interesting, due to examples from \cite{subramani-2023} (see Fig.~\ref{fig:subramani}) showing dynamics different from the $\beta_1=\beta_2$ case analyzed in our current work. This might be an interesting direction to explore in our future research.

\section*{Author contributions (CRediT)}

\textbf{Dorian Falęcki:} Conceptualization, Methodology, Software, Formal analysis, Investigation, Resources, Data Curation, Writing -- Original Draft, Visualization;
\textbf{Mikołaj Rosman:} Conceptualization, Methodology, Software, Formal analysis, Investigation, Resources, Data Curation, Writing -- Original Draft, Writing -- Review \& Editing, Visualization;
\textbf{Michał Palczewski:} Conceptualization, Methodology, Investigation;
\textbf{Paweł Pilarczyk:} Conceptualization, Methodology, Software, Validation, Writing -- Review \& Editing, Visualization, Supervision, Project administration, Funding acquisition;
\textbf{Agnieszka Bartłomiejczyk:} Conceptualization, Methodology, Validation, Writing -- Original Draft, Writing -- Review \& Editing, Supervision, Project administration, Funding acquisition.

\section*{Acknowledgments}

We acknowledge the support obtained from the National Science Centre, Poland, within the grant OPUS 2021/41/B/ST1/00405 (D.~Falęcki, M.~Palczewski, P.~Pilarczyk), as well as the support obtained from Gdańsk University of Technology in the framework of two grants: 9/2020/IDUB/III.4.1/Tc (D.~Falęcki, P.~Pilarczyk) and 27/1/2022/IDUB/III.4c/Tc (M.~Rosman, D.~Falęcki, A.~Bartłomiejczyk), both under the program Technetium Talent Management Grants -- `Excellence Initiative -- Research University'.

The topological--numerical computations discussed in Sec.~\ref{sec:results} were carried out using the computers of Centre of Informatics Tricity Academic Supercomputer \& Network.

\section*{Data availability statement}

The data resulting from the computation discussed in Sec.~\ref{sec:results} is publicly available at the addresses \llabel{dataURLs}\url{https://doi.org/10.34808/86fg-p121} and \url{https://doi.org/10.34808/s5nt-nt30}. An interactive continuation diagram with a browser of Conley--Morse graphs and numerical Morse decompositions is available on \url{https://www.pawelpilarczyk.com/gen2d/}.

\printbibliography

@ARTICLE{Andrecut_2005,
  title     = "Mean field dynamics of random Boolean networks",
  author    = "Andrecut, M",
  journal   = "J. Stat. Mech.",
  publisher = "IOP Publishing",
  volume    =  2005,
  number    =  02,
  pages     = "P02003",
  month     =  feb,
  year      =  2005,
  doi = "10.1088/1742-5468/2005/02/P02003",
}

@ARTICLE{Andrecut_2006,
  title     = "Mean-field model of genetic regulatory networks",
  author    = "Andrecut, M and Kauffman, S A",
  journal   = "New J. Phys.",
  publisher = "IOP Publishing",
  volume    =  8,
  number    =  8,
  pages     = "148--148",
  month     =  aug,
  year      =  2006,
  doi = "10.1088/1367-2630/8/8/148",
}

@ARTICLE{andrecut_main,
  title     = "Chaos in a discrete model of a two-gene system",
  author    = "Andrecut, M and Kauffman, S A",
  journal   = "Phys. Lett. A",
  publisher = "Elsevier BV",
  volume    =  367,
  number    = "4-5",
  pages     = "281--287",
  month     =  07,
  year      =  2007,
  language  = "en",
  doi = "10.1016/j.physleta.2007.03.074",
}

@article{arai-2009,
author = {Z Arai and W Kalies and H Kokubu and K Mischaikow and H Oka and P Pilarczyk},
title = {A Database Schema for the Analysis of Global Dynamics of Multiparameter Systems},
journal = {SIAM Journal on Applied Dynamical Systems},
volume = {8},
number = {3},
pages = {757--789},
year = {2009},
doi = {10.1137/080734935},
}

@article{bartlomiejczyk-bodnar-2020,
  title = {Justification of quasi-stationary approximation in models of gene expression of a self-regulating protein},
  journal = {Commun Nonlinear Sci Numer Simulat},
  volume = {84},
  pages = {105166},
  year = {2020},
  doi = {10.1016/j.cnsns.2020.105166},
  author = {Agnieszka Bartłomiejczyk and Marek Bodnar},
}

@article{bartlomiejczyk-bodnar-2023,
  title = {Hopf bifurcation in time-delayed gene expression model with dimers},
  journal = {Mathematical Methods in the Applied Sciences},
  volume = {46},
  pages = {12087--12111},
  year = {2023},
  doi = {10.1002/mma.8961},
  author = {Agnieszka Bartłomiejczyk and Marek Bodnar},
}

@article{KAUFFMAN1969437,
title = {Metabolic stability and epigenesis in randomly constructed genetic nets},
journal = {Journal of Theoretical Biology},
volume = {22},
number = {3},
pages = {437--467},
year = {1969},
issn = {0022-5193},
doi = {10.1016/0022-5193(69)90015-0},
author = {S A Kauffman},
}

@article{li-cao-li-jin-2018,
  title = {Upregulation of {HES}1 promotes cell proliferation and invasion in breast cancer as a prognosis marker and therapy target via the {AKT} pathway and {EMT} process},
  journal = {Journal of Cancer},
  volume = {9(4)},
  pages = {757--766},
  year = {2018},
  doi = {10.7150/jca.22319},
  author = {Xiaoying Li and Yu Cao and Mu Li and Feng Jin},
}

@article{liu-dai-du-2015,
  title = {{H}es1: a key role in stemness, metastasis and multidrug resistance},
  journal = {Cancer Biology \& Therapy},
  volume = {16(3)},
  pages = {353--359},
  year = {2015},
  doi = {10.1080/15384047.2015.1016662},
  author = {Zi-Hao Liu and Xiao-Meng Dai and Bin Du},
}

@article{luzzatto-2011,
  title     = {Finite resolution dynamics},
  author    = {Luzzatto, Stefano and Pilarczyk, Pawe{\l}},
  journal   = {Found. Comut. Math.},
  publisher = {Springer Science and Business Media LLC},
  volume    =  {11},
  number    =  {2},
  pages     = {211--239},
  month     =  {4},
  year      =  {2011},
  language  = {en},
  doi = {10.1007/s10208-010-9083-z},
}

@article{piotorwska-bartlomiejczyk-bodnar-2018,
  title = {Mathematical analysis of a generalised p53-{M}dm2 protein gene expression model},
  journal = {Applied Mathematics and Computation},
  volume = {328},
  pages = {26--44},
  year = {2018},
  doi = {10.1016/j.amc.2018.01.014},
  author = {Monika J. Piotrowska and Agnieszka Bartłomiejczyk and Marek Bodnar},
}

@article{pilarczyk-2010,
  title = {Parallelization Method for a Continuous Property},
  volume = {10},
  ISSN = {1615-3383},
  doi = {10.1007/s10208-009-9050-8},
  number = {1},
  journal = {Foundations of Computational Mathematics},
  publisher = {Springer Science and Business Media LLC},
  author = {Pilarczyk, Paweł},
  year = {2010},
  pages = {93--114}
}

@article{pilarczyk-2023,
    author = {Pilarczyk, Paweł and Signerska-Rynkowska, Justyna and Graff, Grzegorz},
    title = "{Topological-numerical analysis of a two-dimensional discrete neuron model}",
    journal = {Chaos: An Interdisciplinary Journal of Nonlinear Science},
    volume = {33},
    number = {4},
    pages = {043110},
    year = {2023},
    month = {04},
    issn = {1054-1500},
    doi = {10.1063/5.0129859},
}

@article{pilarczyk-graff-2024,
  title = {An absorbing set for the {C}hialvo map},
  journal = {Communications in Nonlinear Science and Numerical Simulation},
  volume = {132},
  pages = {107947},
  year = {2024},
  issn = {1007--5704},
  doi = {10.1016/j.cnsns.2024.107947},
  author = {Paweł Pilarczyk and Grzegorz Graff},
}

@article{subramani-2023,
author = {Subramani, Rajeshkanna and Kadhim, Hayder and Rajagopal, Karthikeyan and Krejcar, Ondrej and Namazi, Hamidreza},
year = {2023},
month = {07},
pages = {},
title = {The dynamic analysis of discrete fractional-order two-gene map},
journal = {The European Physical Journal Special Topics},
doi = {10.1140/epjs/s11734-023-00912-7},
}

@article{sharma2019,
author = {Sharma, R. and Saha, Lal},
title = {Dynamics of two-gene {A}ndrecut-{K}auffman system: Chaos and complexity},
year = {2019},
month = {01},
journal = {Italian Journal of Pure and Applied Mathematics},
volume = {41},
pages = {405--413},
}

@ARTICLE{De_Souza2012-vs,
  title     = "Self-similarities of periodic structures for a discrete model of a two-gene system",
  author    = "de Souza, S. L. T. and Lima, A. A. and Caldas, I. L. and Medrano-T., R. O. and Guimar{\~a}es-Filho, Z. O.",
  journal   = "Phys. Lett. A",
  publisher = "Elsevier BV",
  volume    =  376,
  number    =  15,
  pages     = "1290--1294",
  month     =  03,
  year      =  2012,
  language  = "en",
  doi = "10.1016/j.physleta.2012.02.036",
}

@article{rosman-2025,
title = {Bistability and chaos in the discrete two-gene Andrecut-Kauffman model},
journal = {Discrete and Continuous Dynamical Systems - B},
pages = {4442-4461},
volume = {30},
number = {11},
year = {2025},
doi = {10.3934/dcdsb.2025028},
author = {Mikołaj Rosman and Michał Palczewski and Paweł Pilarczyk and Agnieszka Bartłomiejczyk}
}

@article{pilarczyk-2008,
title = {Excision-preserving cubical approach to the algorithmic computation of the discrete Conley index},
journal = {Topology and its Applications},
volume = {155},
number = {10},
pages = {1149-1162},
year = {2008},
issn = {0166-8641},
doi = {10.1016/j.topol.2008.02.003},
url = {https://www.sciencedirect.com/science/article/pii/S0166864108000382},
author = {Paweł Pilarczyk and Kinga Stolot},
}

@article{knipl-2015,
author = {Knipl, Di\'{a}na H. and Pilarczyk, Pawe\l{} and R\"{o}st, Gergely},
title = {Rich Bifurcation Structure in a Two-Patch Vaccination Model},
journal = {SIAM Journal on Applied Dynamical Systems},
volume = {14},
number = {2},
pages = {980-1017},
year = {2015},
doi = {10.1137/140993934},
}

@article{miyaji-2016,
title = {A study of rigorous ODE integrators for multi-scale set-oriented computations},
journal = {Applied Numerical Mathematics},
volume = {107},
pages = {34-47},
year = {2016},
issn = {0168-9274},
doi = {10.1016/j.apnum.2016.04.005},
url = {https://www.sciencedirect.com/science/article/pii/S0168927416300435},
author = {Tomoyuki Miyaji and Paweł Pilarczyk and Marcio Gameiro and Hiroshi Kokubu and Konstantin Mischaikow},
}

@InProceedings{GAIO,
author="Dellnitz, Michael
and Froyland, Gary
and Junge, Oliver",
editor="Fiedler, Bernold",
title="The Algorithms Behind GAIO --- Set Oriented Numerical Methods for Dynamical Systems",
booktitle="Ergodic Theory, Analysis, and Efficient Simulation of Dynamical Systems",
year="2001",
publisher="Springer Berlin Heidelberg",
address="Berlin, Heidelberg",
pages="145--174",
isbn="978-3-642-56589-2",
doi="10.1007/978-3-642-56589-2_7"
}

@article{Dellnitz1997,
  title = {A subdivision algorithm for the computation of unstable manifolds and global attractors},
  volume = {75},
  ISSN = {0945-3245},
  url = {http://dx.doi.org/10.1007/s002110050240},
  DOI = {10.1007/s002110050240},
  number = {3},
  journal = {Numerische Mathematik},
  publisher = {Springer Science and Business Media LLC},
  author = {Dellnitz,  Michael and Hohmann,  Andreas},
  year = {1997},
  month = Jan,
  pages = {293--317}
}

@article{Day2004,
author = {Day, S. and Junge, O. and Mischaikow, K.},
title = {A Rigorous Numerical Method for the Global Analysis of Infinite-Dimensional Discrete Dynamical Systems},
journal = {SIAM Journal on Applied Dynamical Systems},
volume = {3},
number = {2},
pages = {117-160},
year = {2004},
doi = {10.1137/030600210}
}

@book{KMM2004,
  title = {Computational Homology},
  ISBN = {9780387215976},
  ISSN = {2196-968X},
  DOI = {10.1007/b97315},
  journal = {Applied Mathematical Sciences},
  publisher = {Springer New York},
  author = {Kaczynski,  Tomasz and Mischaikow,  Konstantin and Mrozek,  Marian},
  year = {2004}
}

@incollection{mischaikow2002conley,
  ISBN = {9780444501684},
  ISSN = {1874-575X},
  url = {http://dx.doi.org/10.1016/S1874-575X(02)80030-3},
  DOI = {10.1016/s1874-575x(02)80030-3},
  booktitle = {Handbook of Dynamical Systems},
  publisher = {Elsevier},
  year = {2002},
  pages = {393--460},
  title     = {Conley index},
  author    = {Mischaikow, Konstantin and Mrozek, Marian},
  volume    = {2},
  editor    = {Fiedler, Bernold},
  isbn      = {978-0-444-50168-4},
  doi       = {10.1016/S1874-575X(02)80030-3}
}

@article{Szymczak1997,
  title = {A combinatorial procedure for finding isolating neighbourhoods and index pairs},
  volume = {127},
  ISSN = {1473-7124},
  DOI = {10.1017/s0308210500026901},
  number = {5},
  journal = {Proceedings of the Royal Society of Edinburgh: Section A Mathematics},
  publisher = {Cambridge University Press (CUP)},
  author = {Szymczak,  Andrzej},
  year = {1997},
  pages = {1075--1088}
}

@book{Moore2009,
  title = {Introduction to Interval Analysis},
  ISBN = {9780898717716},
  DOI = {10.1137/1.9780898717716},
  publisher = {Society for Industrial and Applied Mathematics},
  author = {Moore,  Ramon E. and Kearfott,  R. Baker and Cloud,  Michael J.},
  year = {2009},
  month = Jan 
}

@article{MM1995chaos,
  title = {Chaos in the Lorenz equations: a computer-assisted proof},
  volume = {32},
  ISSN = {0273-0979},
  DOI = {10.1090/s0273-0979-1995-00558-6},
  number = {1},
  journal = {Bulletin of the American Mathematical Society},
  publisher = {American Mathematical Society (AMS)},
  author = {Mischaikow,  Konstantin and Mrozek,  Marian},
  year = {1995},
  pages = {66--72}
}

@article{Mrozek1989,
  title = {Index pairs and the fixed point index for semidynamical systems with discrete time},
  volume = {133},
  ISSN = {1730-6329},
  DOI = {10.4064/fm-133-3-179-194},
  number = {3},
  journal = {Fundamenta Mathematicae},
  publisher = {Institute of Mathematics,  Polish Academy of Sciences},
  author = {Mrozek,  Marian},
  year = {1989},
  pages = {179--194}
}

@article{Mrozek1990,
  title = {Open index pairs,  the fixed point index and rationality of zeta functions},
  volume = {10},
  ISSN = {1469-4417},
  DOI = {10.1017/s0143385700005745},
  number = {3},
  journal = {Ergodic Theory and Dynamical Systems},
  publisher = {Cambridge University Press (CUP)},
  author = {Mrozek,  Marian},
  year = {1990},
  month = Sep,
  pages = {555--564}
}

@article{Kapela2021,
  title = {CAPD::DynSys: A flexible C++ toolbox for rigorous numerical analysis of dynamical systems},
  volume = {101},
  ISSN = {1007-5704},
  DOI = {10.1016/j.cnsns.2020.105578},
  journal = {Communications in Nonlinear Science and Numerical Simulation},
  publisher = {Elsevier BV},
  author = {Kapela,  Tomasz and Mrozek,  Marian and Wilczak,  Daniel and Zgliczyński,  Piotr},
  year = {2021},
  month = Oct,
  pages = {105578}
}

@article{Dey2022,
  title = {Persistence of Conley--Morse Graphs in Combinatorial Dynamical Systems},
  volume = {21},
  ISSN = {1536-0040},
  url = {http://dx.doi.org/10.1137/21M143162X},
  DOI = {10.1137/21m143162x},
  number = {2},
  journal = {SIAM Journal on Applied Dynamical Systems},
  publisher = {Society for Industrial & Applied Mathematics (SIAM)},
  author = {Dey,  Tamal K. and Mrozek,  Marian and Slechta,  Ryan},
  year = {2022},
  month = Apr,
  pages = {817--839}
}

@article{Mischaikow2016,
  title = {Discretization strategies for computing Conley indices and Morse decompositions of flows},
  volume = {3},
  ISSN = {2158-2505},
  url = {http://dx.doi.org/10.3934/jcd.2016001},
  DOI = {10.3934/jcd.2016001},
  number = {1},
  journal = {Journal of Computational Dynamics},
  publisher = {American Institute of Mathematical Sciences (AIMS)},
  author = {Mischaikow,  Konstantin and Mrozek,  Marian and Weilandt,  Frank},
  year = {2016},
  pages = {1--16}
}

@article{Galias2001,
  title = {Interval methods for rigorous investigations of periodic orbits},
  volume = {11},
  ISSN = {1793-6551},
  url = {http://dx.doi.org/10.1142/S0218127401003516},
  DOI = {10.1142/s0218127401003516},
  number = {09},
  journal = {International Journal of Bifurcation and Chaos},
  publisher = {World Scientific Pub Co Pte Lt},
  author = {Galias, Zbigniew},
  year = {2001},
  month = Sep,
  pages = {2427--2450}
}

@article{Galias2021,
  author = {Zbigniew Galias},
  title = {Periodic orbits of the logistic map in single and double precision implementations},
  year = 2021,
  journal = {IEEE Trans. Circuits Syst. II, Exp. Briefs},
  volume = 68,
  number = 11,
  doi = {10.1109/TCSII.2021.3081604},
  pages = {3471--3475},
}

\end{document}